\documentclass[11pt]{article}
 \usepackage{amsfonts}
 \usepackage{mathrsfs}
 \usepackage{bbm}
 \usepackage{amsthm}
\usepackage{amsmath,amssymb}
\def\MatrixFont{\bf}
\def\VectorFont{\bf}

\newcommand{\mA}{{\MatrixFont A}}
\newcommand{\mB}{{\MatrixFont B}}
\newcommand{\mC}{{\MatrixFont C}}
\newcommand{\mD}{{\MatrixFont D}}
\newcommand{\mE}{{\MatrixFont E}}
\newcommand{\mF}{{\MatrixFont F}}

\newcommand{\mI}{{\MatrixFont I}}
\newcommand{\mJ}{{\MatrixFont J}}

\newcommand{\mM}{{\MatrixFont M}}

\newcommand{\mP}{{\MatrixFont P}}
\newcommand{\mQ}{{\MatrixFont Q}}
\newcommand{\mR}{{\MatrixFont R}}

\newcommand{\mX}{{\MatrixFont X}}

\newcommand{\mZ}{{\MatrixFont Z}}

\newcommand{\vh}{{\VectorFont h}}

\newcommand{\vv}{{\VectorFont v}}
\newcommand{\vw}{{\VectorFont w}}
\newcommand{\vx}{{\VectorFont x}}

\newcommand{\vz}{{\VectorFont z}}

\newtheorem{theorem}{Theorem}[section]
\newtheorem{lemma}[theorem]{Lemma}

\theoremstyle{definition}

\theoremstyle{proposition}
\newtheorem{proposition}[theorem]{Proposition}

\theoremstyle{remark}
\newtheorem{remark}[theorem]{Remark}

\theoremstyle{algorithm}

\theoremstyle{corollary}
\newtheorem{corollary}[theorem]{Corollary}

\theoremstyle{example}

\usepackage{amssymb}
\usepackage{epsfig}
\usepackage{times}
\title{Posterior Cram$\acute{e}r$-Rao Bounds for Discrete-Time Nonlinear Filtering with Finitely Correlated Noise}

\author{Zhiguo Wang and Xiaojing Shen
\thanks{This work was supported  in part by the open research funds of BACC-STAFDL of China
under Grant No. 2013afdl011,  the special funds of NEDD of China under Grant No. 201314 and the PCSIRT1273.} 
\thanks{Zhiguo Wang and Xiaojing Shen (corresponding author) are with Department of Mathematics, Sichuan University,
Chengdu, Sichuan 610064, China.  E-mail: shenxj@scu.edu.cn.}}

\begin{document}
 \maketitle
\begin{abstract}
In this paper, an approximation recursive formula of the mean-square error lower bound for the discrete-time nonlinear filtering problem when noises of dynamic systems are temporally correlated is derived based on the Van Trees (posterior) version of the Cram$\acute{e}$r-Rao inequality. The formula is unified in the sense that it can be applied to the multi-step correlated process noise, multi-step correlated measurement noise and multi-step cross-correlated process and measurement noise simultaneously. The lower bound is evaluated by two typical target tracking examples respectively. Both of them show that the new lower bound is significantly different from that of the method which ignores correlation of noises. Thus, when they are applied to sensor selection problems, number of selected sensors becomes very different to obtain a desired estimation performance.
\end{abstract}

\noindent{\bf keywords:} Nonlinear filtering; correlated noises; posterior Cram$\acute{e}$r-Rao bounds; target tracking; sensor networks; sensor
selection

\section{Introduction}\label{sec_1}
The problem of discrete-time nonlinear filtering when noises of dynamic systems are \emph{temporally correlated} (i.e., colored) arises in
various applications such as target tracking, navigation, stochastic approximation, adaptive control, robotics, mobile communication
\cite{Li-Jilkov03,Mazor-Averbuch-BarShalom-Dayan98,Simon06,BarShalom-Li-Kirubarajan01,Ljung-Gunnarsson90,Guo-94,Maryak-Spall-Silberman-95,Rogers87,Halevi90,Blair-Watson-Rice91,Rapoport-Oshman05,Yun-Bachmann06,Jiang-Zhou-Zhu10,Chen12}, just to name a few. For example, in maneuvering
target tracking \cite{Li-Jilkov03}, the process noise and target acceleration can be characterized as temporally correlated stochastic process, respectively. In tracking airborne or missile targets using radar data, the measurement noise is significantly correlated when the measurement
frequency is high \cite{Mazor-Averbuch-BarShalom-Dayan98}. When a system is an airplane and winds are buffeting the plane, an anemometer is used
to measure wind speed as an input to Kalman filter. So the random gusts of wind affect both the process (i.e., the airplane dynamics) and the
measurement (i.e., the sensed wind speed). Thus, there is a correlation between the process noise and the measurement noise \cite{Simon06}. More
detailed results and discussions can be seen in Chapter 7 of the book \cite{Simon06}, Chapter 8 of the book \cite{BarShalom-Li-Kirubarajan01},
and reference therein.
As is well known, the optimal estimator for these problems cannot be obtained for \emph{nonlinear} and \emph{non-Gaussian} dynamic systems in
general. Besides, assessing the achievable performance of suboptimal filtering techniques may be difficult. A main challenge to researchers in
these fields is to find lower bounds corresponding to optimum performance recursively, which give an indication of performance limitations and
can be used to determine whether imposed performance requirements are realistic or not.



The most popular lower bound is the well-known Cram$\acute{e}$r-Rao bound (CRB). In time-invariant statistical models, the estimated parameter
vector is usually considered real-valued (non-random). The lower bound is given by the inverse of the Fisher information matrix (see, e.g, \cite{VanTrees68,BarShalom2014}). Recently,  the authors in \cite{BarShalom2014} discuss the regularity conditions required for the CRB
for real-valued (unknown) parameters to hold. It
is shown that the commonly assumed requirement that the
support of the likelihood function (LF) should be independent
of the parameter to be estimated can be replaced by the much
weaker requirement that the LF is continuous at the end points
of its support. In the time-varying
systems context we deal with here, the estimated
parameter vector is modeled random. A lower bound that is analogous to the CRB for random parameters
was derived in \cite{VanTrees68};
 this bound is also known as the Van Trees version of the CRB, or referred to as posterior CRB (PCRB)
\cite{Tichavsky-Muravchik-Nehorai98}, where the underlying static random system is assumed to satisfy some regularity conditions which are
presented in Section \ref{sec_2}. In addition, a general class of Weiss-Weinstein lower bounds in parameter estimation is derived in \cite{Weinstein-Weiss88} under less restrictive requirement.
The first derivation of a sequential PCRB version applicable to discrete-time dynamic system filtering, the problem addressed in this paper, was done in \cite{Bobrovsky-Zakai75} and then extended in \cite{Taylor79,Galdos80,Chang81,Kerr89}. The most
general form of sequential PCRB for discrete-time nonlinear systems was presented in \cite{Tichavsky-Muravchik-Nehorai98}, \cite{Koshaev-Stepanov97}. Together with the
original static form of the CRB, these results served as a basis for a large number of applications
\cite{Schultheiss-Weinstein81,Aidala-Hammel83,Kirubarajan-BarShalom96,
Niu-Willett-BarShalom01,Zhang-Willett-BarShalom02,Zheng-Ozdemir-Niu-Varshney12,Kar-Varshney-Palaniswami12}.

In this paper, we focus on an approximation recursive derivation of the PCRB for the discrete-time nonlinear and non-Gaussian filtering problem when noises of
dynamic systems are temporally correlated. The derived formula is unified in the sense that it can be applied to the multi-step correlated
process noise, multi-step correlated measurement noise and multi-step cross-correlated process and measurement noise simultaneously.
The derivation differs from the other approaches that instead consider the three cases separately and assume the linear or Gaussian dynamic
systems. Although the unified formula can come across the three cases of finite-step correlated noises, a few corollaries with simpler formulae
follow to elucidate special cases, which may be used more frequently. The main results are presented in Section \ref{sec_3}.
In Section IV, the new lower bound is evaluated by two typical target tracking examples, respectively. Both of them show that the new lower bound
is significantly different from that of the existing methods. Thus, when they are applied to sensor selection
problems, simulations show that the new formula can derive a more accurate number of selected sensors to obtain a desired estimation performance.
Conclusions are drawn in Section \ref{sec_5}. In order to enhance readability, all proofs are given in Appendices.

\section{Preliminaries}\label{sec_2}
\subsection{Problem Formulation}\label{sec_2_1}
Consider a nonlinear dynamic system 
\begin{eqnarray}
\label{Eqs_sec2_1} \vx_{k+1}&=&f_k(\vx_k, \vw_k)\\
\label{Eqs_sec2_2} \vz_k&=&h_k(\vx_k, \vv_k)
\end{eqnarray}
where $\vx_{k}\in \mathbb{R}^r$ is the state to be estimated at time $k$, $r$ is the dimension of the state; $\vz_k\in \mathbb{R}^n$ is the
measurement vector. The function $f_k$ and $h_k$ are nonlinear functions in general. $\{\vw_k\}$ and $\{\vv_k\}$ are noises both temporally
\emph{finite-step correlated}, respectively. We discuss the following three cases:
\begin{enumerate}
\item The process noises are \emph{$l$-step auto-correlated} if their probability density functions satisfy
$p(\vw_k$, $\vw_{k-i})\neq p(\vw_k)$ $p(\vw_{k-i})$, $p(\vw_k,\vw_{k-j})=p(\vw_k)p(\vw_{k-j})$, for ~$i=1,\ldots,l$, $j=l+1,\ldots, k$, $k\geq
l+1$. We denote by $0$-step correlated process noise if they are temporally independent.

\item The measurement noises are \emph{$l$-step auto-correlated} if
$p(\vv_k,\vv_{k-i})\neq p(\vv_k)p(\vv_{k-i})$, and $p(\vv_k,\vv_{k-j})$ $=p(\vv_k)p(\vv_{k-j})$, for ~$i=1,\ldots,l$, $j=l+1,\ldots, k$, $k\geq
l+1$. We denote by $0$-step correlated measurement noise if they are temporally independent.

\item The measurement noise is \emph{backward $l$-step cross-correlated} with
process noise, if $p(\vv_k$, $\vw_{k-i})\neq p(\vv_k)p(\vw_{k-i})$, and $p(\vv_k,\vw_{k-j})=p(\vv_k)p(\vw_{k-j})$, for \,$i=1,\ldots,l,
j=l+1,\ldots, k$, $k\geq l+1$; The measurement noise is \emph{forward $l$-step cross-correlated} with process noise, if
$p(\vv_k,\vw_{k-1+i})\neq p(\vv_k)$ $p(\vw_{k-1+i})$, and $p(\vv_k,\vw_{k-1+j})=p(\vv_k)p(\vw_{k-1+j})$, for $i=1,\ldots,l, j=l+1,\ldots, k$,
$k\geq l+1$; We denote that the measurement noise is forward and backward $0$-step correlated with the process noise if they are mutually
independent.
\end{enumerate}
Note that if the measurement noise is finite-step correlated to the process noise, then the process noise is also finite-step correlated to the
measurement noise. For example, the observation noises at times $k-1$ and $k$ are correlated to process noise at
time $k-2$, then process noises at times $k-2$ and $k-1$ are correlated to observation noise at time $k$. Thus, their recursive formulae are the same and we only consider the former.

Since, in target tracking, the three correlated cases may be encountered simultaneously ( see, e.g.,
\cite{Li-Jilkov03,Mazor-Averbuch-BarShalom-Dayan98}) and the optimal estimator for these problems cannot be obtained for \emph{nonlinear} and
\emph{non-Gaussian} dynamic systems in general, the
goal of this paper is to derive a unified lower bound recursively, which can be
used to determine whether imposed performance requirements
are realistic or not.

\subsection{Posterior Cram$\acute{e}$r-Rao Bounds}\label{sec_2_2}

Let $\vx$ be a $r$-dimensional random parameter vector and $\vz$ be a measurement vector, let $p_{\vx,\vz}(\mX,\mZ)$ be a joint density of the pair
$(\vx,\vz)$. The mean-square error of any estimate $\hat{\vx}(\mZ)$ of $\vx$ satisfies the inequality
\begin{eqnarray}
\label{Eqs_sec2_3} \mP\triangleq E\{[\hat{\vx}(\mZ)-\vx][\hat{\vx}(\mZ)-\vx]^T\}\geq \mJ^{-1}
\end{eqnarray}
where $\mJ$ is the $r\times r$ (Fisher) information matrix with the elements
\begin{eqnarray}
\label{Eqs_sec2_4}\mJ_{ij}=E[-\frac{\partial^2 \ln~p_{\vx,\vz}(\mX,\mZ)}{\partial \mX_i\partial \mX_j}]~~~~i, j=1, \ldots, r
\end{eqnarray}
and the expectation is over both $\vx$ and $\vz$. The superscript $``^T"$ in (\ref{Eqs_sec2_3}) denotes the transpose of a matrix. The following
conditions are assumed to exist:

\begin{enumerate}
\item $\frac{\partial \ln~p_{\vx,\vz}(\mX,\mZ)}{\partial \mX}$ is absolutely integrable with respect to $\mX$ and $\mZ$.
\item $\frac{\partial^2 \ln~p_{\vx,\vz}(\mX,\mZ)}{\partial \mX^2}$ is absolutely integrable with respect to $\mX$ and $\mZ$.
\item The conditional expectation of the error, given $\mX$, is
\begin{eqnarray}
\nonumber B(\mX)=\int_{-\infty}^{+\infty}[\hat{\vx}(\mZ)-\mX]p_{\vz|\vx}(\mZ|\mX)d\mZ
\end{eqnarray}
and assume that
\begin{eqnarray}
\nonumber\lim_{\mX_i\rightarrow\infty}B(\mX)p_\vx(\mX)&=&0, \lim_{\mX_i\rightarrow-\infty}B(\mX)p_\vx(\mX)=0,\\ \nonumber for ~i&=&1,\ldots,r.
\end{eqnarray}
\end{enumerate}
The proof is given in \cite{VanTrees68}.

Assume now that the parameter $\vx$ is decomposed into two parts as $\vx=[\vx_1^T,\vx_2^T]^T$, and the information matrix $\mJ$ is
correspondingly decomposed into blocks

\begin{eqnarray}
\label{Eqs_sec2_5} \mJ=\left(
\begin{array}{cc}
\mJ_{11} &\mJ_{12}\\
\mJ_{12}^T & \mJ_{22} \\
\end{array}
\right)
\end{eqnarray}
It can easily be shown that the covariance of estimation of $\vx_2$ is lower bounded by the right-lower block of $\mJ^{-1}$, i.e.,
\begin{eqnarray}\label{Eqs_sec2_6}
\nonumber \mP_2&\triangleq&
E\{[\hat{\vx}_2(\mZ)-\vx_2][\hat{\vx}_2(\mZ)-\vx_2]^T\}\\[3mm]
&\geq&[\mJ_{22}-\mJ_{12}^T\mJ_{11}^{-1}\mJ_{12}]^{-1}
\end{eqnarray}
assuming that $\mJ_{11}^{-1}$ exists. The matrix $\mJ_{22}-\mJ_{12}^T\mJ_{11}^{-1}\mJ_{12}$ is called the information submatrix for parameter
$\vx_2$ in \cite{Tichavsky-Muravchik-Nehorai98}.

Note that
the joint probability densities
$p(\mX_k, \mZ_k)$ of $\mX_k=(\vx_0^T, \ldots, \vx_k^T)^T$ and $\mZ_k=(\vz_0^T, \ldots, \vz_k^T)^T$ for an arbitrary $k$ is is determined by Equations (\ref{Eqs_sec2_1}) and (\ref{Eqs_sec2_2}) together with $p(\vx_0)$ and also by the noise pdfs
\cite{Bobrovsky-Zakai75}. The conditional probability densities $p(\vx_{k+1}|\mX_k,\mZ_k)$ and $p(\vz_{k+1}|\mX_{k+1},\mZ_k)$ can be obtained
from (\ref{Eqs_sec2_1}) and (\ref{Eqs_sec2_2}), respectively, under suitable hypotheses. In this paper, $p(\mX_k,\mZ_k)$ is denoted by $p_k$ for brevity. From a Bayesian perspective, the joint probability function of $\mX_{k+1}$ and $\mZ_{k+1}$ can be written as
\begin{eqnarray}
\label{Eqs_sec2_7}p_{k+1}&=&p(\mX_{k+1},\mZ_{k+1})\\
\label{Eqs_sec2_8}&=&p_kp(\vx_{k+1}|\mX_k,\mZ_k)p(\vz_{k+1}|\mX_{k+1},\mZ_k).
\end{eqnarray}
In addition, define $\nabla$ and $\triangle$ be the first and second-order operator partial derivatives, respectively
\begin{eqnarray}
\label{Eqs_sec2_07} \nabla_\alpha&=&\left( \begin{array}{cccc} \frac{\partial}{\partial\alpha_1}
& \frac{\partial}{\partial\alpha_2} & \ldots & \frac{\partial}{\partial\alpha_n} \\
\end{array}
\right)^T~~~~\forall\alpha\in R^n,\\
\label{Eqs_sec2_08} \triangle_\alpha^{\beta}&=&\nabla_\alpha\nabla_\beta^T.
\end{eqnarray}
Using this notation and (\ref{Eqs_sec2_7}), (\ref{Eqs_sec2_4}) can be written as
\begin{eqnarray}
\label{Eqs_sec2_9} \mJ(\mX_k)=E\left(-\triangle_{X_k}^{X_k}\ln p(\mX_k,\mZ_k)\right)=E\left(-\triangle_{X_k}^{X_k}\ln ~p_k\right).
\end{eqnarray}
Decompose state vector $\mX_k$ as $\mX_k$=$(\mX_{k-1}^T,\vx_k^T)^T$ and the ($kr\times kr$) information matrix $J(\mX_k)$ correspondingly as
\begin{eqnarray}
\label{Eqs_sec2_10} \mJ(\mX_k)=\left(
\begin{array}{cc}
\mA_k^{11} &\mA_k^{12}\\
{\mA_k^{21}} & \mA_k^{22} \\
\end{array}
\right)
\end{eqnarray}
where
\begin{eqnarray}\label{Eqs_sec2_11}
\nonumber\mA_k^{11}&=&E(-\bigtriangleup_{\mX_{k-1}}^{\mX_{k-1}}\ln~p_k),\\
 \nonumber\mA_k^{12}&=&E(-\bigtriangleup_{\mX_{k-1}}^{\vx_k}\ln~p_k),\\
  \nonumber\mA_k^{21}&=&E(-\bigtriangleup_{\vx_k}^{\mX_{k-1}}\ln~p_k),\\
\nonumber\mA_k^{22}&=&E(-\bigtriangleup_{\vx_k}^{\vx_k}\ln~p_k).
\end{eqnarray}
Thus, the posterior information submatrix for estimating $\vx_k$, denoted by $\mJ_k$, which is given as the inverse of the ($r\times r$) right
block of $[\mJ(\mX_k)]^{-1}$, i.e.,
\begin{eqnarray}
\label{Eqs_sec2_12}\mJ_k=\mA_k^{22}-{\mA_k^{21}}{\mA_k^{11}}^{-1}\mA_k^{12}.
\end{eqnarray}
$\mJ_k^{-1}$ is the PCRB of estimating state vector $\vx_k$.

In the following, we derive the recursive formula of the posterior information submatrices $\{\mJ_k\}$ when the noises of dynamic systems are
finite-step correlated.

\section{Main results}\label{sec_3}
In this section, we address the recursive formula of the posterior information submatrices $\{\mJ_k\}$ when the noises of dynamic systems are
finite-step correlated.
Let us give some remarks about the general matrix $M_k$ on the notation: 
\begin{enumerate}
\item $\mM_k^{i,j}$, $i,j=1,\ldots,l$ denotes the $i$-th row and $j$-th column block
of the $(l\times l)$ block matrix $\mM_k$ at time $k$;
\item If $i\leq 0$, or $j\leq 0$, then $\mM_k^{i,j}=0$;
\item If $i>l$, or $j>l$, then $\mM_k^{i,j}=0$.
\end{enumerate}
When the measurement noise and the process noise of the dynamic system (\ref{Eqs_sec2_1})-(\ref{Eqs_sec2_2}) are  temporally auto-correlated and
cross-correlated simultaneously, we have the following unified recursion as follows.

\begin{proposition}\label{thm_1}
If the measurement noise is $l_1$-step auto-correlated ($l_1\geq0$), the process noise is $l_2$-step auto-correlated ($l_2\geq0$), and the measurement
noise is backward $l_3$-step and forward $l_4$-step cross-correlated with the process noise ($l_3\geq0$, $l_4\geq0$), then the sequence
$\{\mJ_k\}$ of posterior information submatrices for estimating state vector $\{\vx_k\}$ approximately obeys the recursion
\begin{eqnarray}
\label{Eqs_sec3_1}\mJ_{k+1}=\mD_k^{22}-\mD_k^{21}(\mD_k^{11}+\mE_k)^{-1}\mD_k^{12},
\end{eqnarray}
where the recursive terms $\mE_k$, $\mD_k^{11}$, $\mD_k^{12}$, $\mD_k^{21}$ and $\mD_k^{22}$ are calculated as the following three cases.
In order to facilitate the discussion, with a slight abuse of notations, we denote by $l_3^\prime\triangleq \max\{l_3,1\},l_2^\prime\triangleq
\max\{l_2,1\}$.
\begin{enumerate}
\item  $l_3^\prime>l_2^\prime+1$:

The $i$-th row and $j$-th column block of the matrix $\mE_k$ are recursively calculated as
\begin{eqnarray}
\nonumber \mE_k^{i,j}&=&\mE_{k-1}^{i+1,j+1}+\mB_{k-1}^{i+l_2^\prime-l_3^\prime+2,j+l_2^\prime-l_3^\prime+2}\\[3mm]
\label{Eqs_sec3_2}&&+\mC_{k-1}^{i+1,j+1}-(\mE_{k-1}^{i+1,1}+\mC_{k-1}^{i+1,1})\\[3mm]
\nonumber&&\cdot(\mE_{k-1}^{1,1}+\mC_{k-1}^{1,1})^{-1}(\mE_{k-1}^{1,j+1}+\mC_{k-1}^{1,j+1})
\end{eqnarray}
\begin{center}
$for ~i=1,2,\ldots,l_3^\prime-1,~ j=1,2,\ldots,l_3^\prime-1$
\end{center}
where
\begin{eqnarray}
\nonumber&&\mB_{k}^{i,j}=\\
\nonumber&& \left\{
\begin{array}{ll}
E(-\bigtriangleup_{\vx_{k-l_2^\prime+i}}^{\vx_{k-l_2^\prime+j}}\ln~p(\vx_{k+1}
|\vx_k,\ldots,\vx_{k-l_2^\prime+1})) & ~ \\
~~~for ~i,j=1,2,\ldots,l_2^\prime+1,~l_4=0,&~\\[3mm]
E(-\bigtriangleup_{\vx_{k-l_2^\prime+i}}^{\vx_{k-l_2^\prime+j}}\ln~p(\vx_{k+1}|
\vx_k,\ldots,\vx_{k-l_2^\prime+1},&~\\
\qquad\qquad\qquad\qquad\qquad \ \vz_k,\ldots,\vz_{k-l_4+1}))&~\\
~~~for~i,j=1,2,\ldots,l_2^\prime+1,~l_4>0,&~\\
\end{array}
\right.
\end{eqnarray}
\begin{eqnarray}
\label{Eqs_sec3_3}
\end{eqnarray}
\begin{eqnarray}
\nonumber&&\mC_{k}^{i,j}=\\
\nonumber&&\left\{
\begin{array}{ll}
E(-\bigtriangleup_{\vx_{k+1+i-l_3^\prime}}^{\vx_{k+1+j-l_3^\prime}}\ln~p(\vz_{k+1}|\vx_{k+1},\ldots,\vx_{k-l_3^\prime+2})) & ~ \\
~~~for ~i,j=1,2\ldots,l_3^\prime,~l_1=0&~\\[3mm]
E(-\bigtriangleup_{\vx_{k+1+i-l_3^\prime}}^{\vx_{k+1+j-l_3^\prime}}\ln~p(\vz_{k+1}|\vx_{k+1},\ldots,\vx_{k-l_3^\prime+2},&~\\
\qquad\qquad\qquad\qquad\qquad\qquad \ \vz_k,\ldots,\vz_{k-l_1+1}))&~\\
~~~for~i,j=1,2\ldots,l_3^\prime,~l_1>0.&~\\
\end{array}
\right.
\end{eqnarray}
\begin{eqnarray}
\label{Eqs_sec3_4}
\end{eqnarray}
$\mD_k^{11}$, $\mD_k^{12}$, $\mD_k^{21}$ and $\mD_k^{22}$ in (\ref{Eqs_sec3_1}) can be calculated as follows:

\begin{eqnarray}
\label{Eqs_sec3_5}\mD_k^{22}&=&\mC_k^{l_3^\prime,l_3^\prime}+\mB_k^{l_2^\prime+1,l_2^\prime+1},\qquad\qquad\qquad\qquad\qquad
\end{eqnarray}
\begin{eqnarray}
\label{Eqs_sec3_6}\mD_k^{21}=\left(
\begin{array}{ccccc}
\mC_k^{l_3^\prime1} & \cdots & \mC_k^{l_3^\prime,l_3^\prime-l_2^\prime-1}&\cdots&\mC_k^{l_3^\prime,l_3^\prime-1}\\
~&~&+\mB_k^{l_2^\prime+1,1}&~&+\mB_k^{l_2^\prime+1,l_2^\prime}
\end{array}
\right)=(\mD_k^{12})^T,
\end{eqnarray}
\begin{eqnarray}
\label{Eqs_sec3_7}\mD_k^{11}=\left( \begin{array}{ccccc}
\mC_k^{1,1} &\cdots & \mC_k^{1,l_3^\prime-l_2^\prime}&\cdots&\mC_k^{1,l_3^\prime-1} \\
\vdots &\cdots& \vdots &\cdots& \vdots \\
\mC_k^{l_3^\prime-l_2^\prime,1}&\cdots&\mC_k^{l_3^\prime-l_2^\prime,l_3^\prime-l_2^\prime}&\cdots&\mC_k^{l_3^\prime-l_2^\prime,l_3^\prime-1}\\
~&~&+\mB_k^{1,1}&~&+\mB_k^{1,l_2^\prime}\\
\vdots &\cdots& \vdots &\cdots& \vdots \\
\mC_k^{l_3^\prime-1,1} & \cdots & \mC_k^{l_3^\prime-1,l_3^\prime-l_2^\prime}&\cdots&\mC_k^{l_3^\prime-1,l_3^\prime-1} \\
~&~&+\mB_k^{l_2^\prime,1}&~&+\mB_k^{l_2^\prime,l_2^\prime}
\end{array}
\right).
\end{eqnarray}
%
\item $l_3^\prime< l_2^\prime+1$

The $i$-th row and $j$-th column block of the matrix $\mE_k$ are recursively calculated as
\begin{eqnarray}
\nonumber \mE_k^{i,j}&=&\mE_{k-1}^{i+1,j+1}+\mC_{k-1}^{i+l_3^\prime-l_2^\prime,j+l_3^\prime-l_2^\prime}+\mB_{k-1}^{i+1,j+1}\\[3mm]
\nonumber&&-(\mE_{k-1}^{i+1,1}+\mB_{k-1}^{i+1,1})(\mE_{k-1}^{1,1}+\mB_{k-1}^{1,1})^{-1}\\[3mm]
&&\cdot(\mE_{k-1}^{1,j+1}+\mB_{k-1}^{1,j+1})\\[3mm]
\nonumber&&for ~i=1,2,\ldots,l_2^\prime, ~ j=1,2,\ldots,l_2^\prime,
\end{eqnarray}
where $\mC_{k}^{i,j}$ and $\mB_{k}^{i,j}$ are defined in (\ref{Eqs_sec3_3})-(\ref{Eqs_sec3_4}). $\mD_k^{11}$, $\mD_k^{12}$, $\mD_k^{21}$ and
$\mD_k^{22}$ in (\ref{Eqs_sec3_1}) can be calculated as follows:

\begin{eqnarray}
\label{Eqs_sec3_8}\mD_k^{22}=\mC_k^{l_3^\prime,l_3^\prime}+\mB_k^{l_2^\prime+1,l_2^\prime+1},\qquad\qquad\qquad~~~~~~
\end{eqnarray}
\begin{eqnarray}
\mD_k^{21}= \left(
\begin{array}{ccccc}
\mB_k^{l_2^\prime+1,1} & \cdots & \mB_k^{l_2^\prime+1,l_2^\prime-l_3^\prime}&\cdots&
\mB_k^{l_2^\prime+1,l_2^\prime}\\
~&~&+\mC_k^{l_3^\prime,1}&~&+\mC_k^{l_3^\prime,l_3^\prime-1}
\end{array}
\right)=(\mD_k^{12})^T,
\end{eqnarray}
\begin{eqnarray}\label{Eqs_sec3_10}
\mD_k^{11}=\left(
  \begin{array}{ccc}
    \mB_k^{1,1} & \ldots & \mB_k^{1,l_2^\prime} \\
    \vdots & \ddots & \vdots \\
    \mB_k^{l_2^\prime,1}  & \ldots & \mB_k^{l_2^\prime,l_2^\prime}+\mC_k^{l_3^\prime-1,1_3^\prime-1} \\
  \end{array}
\right).
\end{eqnarray}\\

\item $l_3^\prime= l_2^\prime+1$

The $i$-th row and $j$-th column block of the matrix $\mE_k$ are recursively calculated as
\begin{eqnarray}
\nonumber \mE_k^{i,j}&=&\mE_{k-1}^{i+1,j+1}+\mC_{k-1}^{i+1,j+1}+\mB_{k-1}^{i+1,j+1}\\[3mm]
\nonumber&&-(\mE_{k-1}^{i+1,1}+\mB_{k-1}^{i+1,1}+\mC_{k-1}^{i+1,1})\\[3mm]
\nonumber&&\cdot(\mE_{k-1}^{1,1}+\mB_{k-1}^{1,1}+\mC_{k-1}^{1,1})^{-1}\\[3mm]
\label{Eqs_sec3_11}&&\cdot(\mE_{k-1}^{1,j+1}+\mB_{k-1}^{1,j+1}+\mC_{k-1}^{1,j+1})
\end{eqnarray}
\begin{center}
$for ~i=1,2,\ldots,l_3^\prime-1, ~ j=1,2,\ldots,l_3^\prime-1$,
\end{center}

where $\mC_{k}^{i,j}$ and $\mB_{k}^{i,j}$ are defined in(\ref{Eqs_sec3_3})-(\ref{Eqs_sec3_4}). $\mD_k^{11}$, $\mD_k^{12}$, $\mD_k^{21}$ and
$\mD_k^{22}$ in (\ref{Eqs_sec3_1}) can be calculated as follows:
\begin{eqnarray}
\label{Eqs_sec3_12}\mD_k^{22}&=&\mC_k^{l_3^\prime,l_3^\prime}+\mB_k^{l_3^\prime,l_3^\prime},\\
\label{Eqs_sec3_13}\nonumber\mD_k^{21}&=&\left(
\begin{array}{ccc}
\mB_k^{l_3^\prime,1} & \cdots &
\mB_k^{l_3^\prime,l_3^\prime-1}\\
+\mC_k^{l_3^\prime,1}&~&+\mC_k^{l_3^\prime,l_3^\prime-1}\\
\end{array}
\right)\\
&=&(\mD_k^{12})^T,\\
\label{Eqs_sec3_14}\mD_k^{11}&=&\left( \begin{array}{ccccc}
\mB_k^{1,1} &\cdots&\mB_k^{1,l_3^\prime-1} \\
+ \mC_k^{1,1}&~&+\mC_k^{1,l_3^\prime-1}\\
\cdots & \ddots &\cdots \\
\mB_k^{l_3^\prime-1,1} & \cdots & \mB_k^{l_3^\prime-1,l_3^\prime-1}\\
+\mC_k^{l_3^\prime-1,1}&~&+\mC_k^{l_3^\prime-1,l_3^\prime-1} \\
\end{array}
\right).
\end{eqnarray}

\end{enumerate}
\end{proposition}
Proof: See the Appendix.

\begin{remark}
The difficulty of the recursion is the derivation of the recursive matrix $\mE_k$, which thanks to two lemmas about the inverse of a matrix
given in Appendix and Schur complement. Although the derivation is very complicated, the final formula is not complicated. The finite-step
correlation of noises is used to determine (\ref{Eqs_sec3_3}) and (\ref{Eqs_sec3_4}), which can be calculated  by analytical or numerical methods.
Note that both of them are approximation equations (see proof in in the appendix).
The initial information submatrix $\mE_0^{ij}$ can be calculated from the a priori probability function $p(\mX_{l_{max}},\mZ_{l_{max}})$ where
$l_{max}=\max\{l_1,l_2,l_3,l_4\}$.
\end{remark}

Although the unified formula can come across the three cases of finite-step correlated noises, a few corollaries with simpler formulae follow to
elucidate special cases, which may be used more frequently.

\begin{corollary}\label{coro_2}
If the measurement noise is backward $l$-step cross-correlated with the process noise ($l\geq0$), the process noise and measurement noise are
temporally independent, respectively, i.e., $l_1=0$, $l_2=0$, $l_3=l$, $l_4=0$, then the sequence $\{\mJ_k\}$ of posterior information
submatrices for estimating state vector $\{\vx_k\}$ obeys the recursion
\begin{eqnarray}
\label{Eqs_sec3_15}\mJ_{k+1}=\mD_k^{22}-\mD_k^{21}(\mD_k^{11}+\mE_k)^{-1}\mD_k^{12},
\end{eqnarray}
where the recursive term $\mE_k$ is calculated as follows:
\begin{enumerate}
\item If $l=1$ or $l=0$, then
\begin{eqnarray}
\nonumber\mE_k&=&
\mB_{k-1}^{2,2}+\mC_{k-1}^{1,1}\\[3mm]
\label{Eqs_sec3_16}&&-\mB_{k-1}^{2,1}(\mE_{k-1}+\mB_{k-1}^{1,1})^{-1}\mB_{k-1}^{1,2}
\end{eqnarray}
where
\begin{eqnarray}
\nonumber \mB_{k-1}^{1,1}&=& E(-\bigtriangleup_{\vx_{k-1}}^{\vx_{k-1}}\ln~p(\vx_{k}|\vx_{k-1})),\\
\nonumber \mB_{k-1}^{1,2}&=&E(-\bigtriangleup_{\vx_{k-1}}^{\vx_{k}}\ln~p(\vx_{k}|\vx_{k-1})),\\
\nonumber \mB_{k-1}^{2,1}&=& E(-\bigtriangleup_{\vx_{k}}^{\vx_{k-1}}\ln~p(\vx_{k}|\vx_{k-1})),\\
 \nonumber\mC_{k-1}^{1,1}&=&E(-\bigtriangleup_{\vx_{k}}^{\vx_{k}}\ln~p(\vz_{k}|\vx_{k})).
\end{eqnarray}
$\mD_k^{11}$, $\mD_k^{12}$, $\mD_k^{21}$ and $\mD_k^{22}$ in (\ref{Eqs_sec3_15}) are  calculated as follows
\begin{eqnarray}
\nonumber \mD_k^{22}&=&\mB_k^{2,2}+\mC_k^{1,1},\\
\nonumber \mD_k^{21}&=&\mB_k^{2,1}=(\mD_k^{12})^T,\\
\nonumber\mD_k^{11}&=&\mB_k^{1,1}.
\end{eqnarray}

\item If $l=2$, then

\begin{eqnarray}
\label{Eqs_sec3_17} \mE_k&=&\mC_{k-1}^{2,2}+\mB_{k-1}^{2,2}-(\mB_{k-1}^{2,1}+\mC_{k-1}^{2,1})\\[3mm]
\nonumber&&\cdot(\mE_{k-1}+\mB_{k-1}^{1,1}+\mC_{k-1}^{1,1})^{-1}(\mC_{k-1}^{1,2}+\mB_{k-1}^{1,2}),
\end{eqnarray}
where
\begin{eqnarray}
\nonumber \mB_{k-1}^{m,n}&=& E(-\bigtriangleup_{\vx_{k-2+m}}^{\vx_{k-2+n}}\ln~p(\vx_{k}|\vx_{k-1}))\\
\nonumber~~&&for~~ m,n=1,2,\\
\nonumber \mC_{k-1}^{m,n}&=& E(-\bigtriangleup_{\vx_{k-2+m}}^{\vx_{k-2+n}}\ln~p(\vz_{k}|\vx_{k},\vx_{k-1}))\\
\nonumber~~&&for~~ m,n=1,2.
\end{eqnarray}
$\mD_k^{11}$, $\mD_k^{12}$, $\mD_k^{21}$ and $\mD_k^{22}$ in (\ref{Eqs_sec3_15}) are  calculated as follows
\begin{eqnarray}
\nonumber\mD_k^{22}&=&\mC_k^{2,2}+\mB_k^{2,2},\\
\nonumber\mD_k^{21}&=&\mB_k^{2,1}+\mC_k^{2,1}=(\mD_k^{12})^T,\\
\nonumber\mD_k^{11}&=&\mB_k^{1,1}+\mC_k^{1,1}.
\end{eqnarray}

\item If $l>2$, then the $i$-th row and $j$-th column block of the matrix $\mE_k$ is recursively calculated as

\begin{eqnarray}
\nonumber\mE_k^{i,j}&=&\mE_{k-1}^{i+1,j+1}+\mB_{k-1}^{i+3-l,j+3-l}+\mC_{k-1}^{i+1,j+1}-\\
\label{Eqs_sec3_18}&&(\mE_{k-1}^{i+1,1}+\mC_{k-1}^{i+1,1})(\mE_{k-1}^{1,1}+\mC_{k-1}^{1,1})^{-1}\\
\nonumber&&\cdot(\mE_{k-1}^{1,j+1}+\mC_{k-1}^{1,j+1}),
\end{eqnarray}
\begin{center}
$for ~i=1,2,\ldots,l_3'-1,~~ j=1,2,\ldots,l_3'-1$
\end{center}
where
\begin{eqnarray}
\nonumber \mC_{k-1}^{m,n}&=& E(-\bigtriangleup_{\vx_{k+m-l}}^{\vx_{k+n-l}}\ln~p(\vz_{k}|\vx_{k},\ldots,\vx_{k+1-l}))\\
\nonumber~~&&for~ m,n=1,2\ldots,l.
\end{eqnarray}
\begin{eqnarray}
\nonumber\mB_{k-1}^{1,1}&=&E(-\bigtriangleup_{\vx_{k-1}}^{\vx_{k-1}}\ln~p(\vx_{k}|\vx_{k-1})),\\
\nonumber\mB_{k-1}^{1,2}&=&E(-\bigtriangleup_{\vx_{k-1}}^{\vx_{k}}\ln~p(\vx_{k}|\vx_{k-1}))=(\mB_k^{2,1})^T,\\
\nonumber\mB_{k-1}^{2,2}&=&E(-\bigtriangleup_{\vx_{k}}^{\vx_{k}}\ln~p(\vx_{k}|\vx_{k-1})).
\end{eqnarray}
$\mD_k^{11}$, $\mD_k^{12}$, $\mD_k^{21}$ and $\mD_k^{22}$ in (\ref{Eqs_sec3_15}) are  calculated as follows
\begin{eqnarray}
\nonumber\mD_k^{22}&=&\mC_k^{l,l}+\mB_k^{2,2},\\
\nonumber\mD_k^{21}&=&\left(
\begin{array}{ccc}
\mC_k^{l,1} & \cdots &
\mB_k^{2,1}+\mC_k^{l,l-1}\\
\end{array}
\right)=(\mD_k^{12})^T,\\
\nonumber\mD_k^{11}&=&\left( \begin{array}{ccc}
\mC_k^{1,1} &\cdots&\mC_k^{1,l-1} \\
\cdots & \ddots& \cdots \\
\mC_k^{l-1,1} & \cdots & \mB_k^{1,1}+\mC_k^{l-1,l-1} \\
\end{array}
\right).
\end{eqnarray}
\end{enumerate}
\end{corollary}
Proof: See the Appendix.

\begin{corollary}\label{coro_3}
If the process noise is $l$-step auto-correlated ($l\geq0$), the measurement noise is temporally independent, and the process noise and the
measurement noise are mutually independent, i.e., $l_1=0$, $l_2=l$, $l_3=0$, $l_4=0$, then the sequence $\{\mJ_k\}$ of posterior information
submatrices for estimating state vector $\{\vx_k\}$ obeys the recursion
\begin{eqnarray}
\label{Eqs_sec3_19} \mJ_{k+1}&=&\mD_k^{22}-\mD_k^{21}(\mE_k+\mD_k^{11})^{-1}\mD_k^{12}
\end{eqnarray}
where the $i$-th row and $j$-th column block of the matrix $\mE_k$ is calculated as follows:
\begin{eqnarray}
\label{Eqs_sec3_20}
\nonumber\mE_k^{i,j}&=&\mE_{k-1}^{i+1,j+1}+\mC_{k-1}^{i+1-l_2^\prime,j+1-l_2^\prime}+\mB_{k-1}^{i+1,j+1}\\[3mm]
&&-(\mE_{k-1}^{i+1,1}+\mB_{k-1}^{i+1,1})(\mE_{k-1}^{1,1}+\mB_{k-1}^{1,1})^{-1}\\[3mm]
\nonumber&&\cdot(\mE_{k-1}^{1,j+1}+\mB_{k-1}^{1,j+1})
\end{eqnarray}
\begin{center}
$for~i=1,2,\ldots,l_2^\prime,~j=1,2,\ldots,l_2^\prime$
\end{center}
where
\begin{eqnarray}
\label{Eqs_sec3_019}\mB_{k-1}^{i,j}=E(-\bigtriangleup_{\vx_{k+i-l_2^\prime-1}}^{\vx_{k+j-l_2^\prime-1}}\ln~p(\vx_{k}|\vx_{k-1},\ldots,\vx_{k-l_2^\prime}))
\end{eqnarray}
\begin{eqnarray}
\label{Eqs_sec3_020}\mC_{k-1}^{1,1}&=&E(-\bigtriangleup_{\vx_{k}}^{\vx_{k}}\ln~p(\vz_{k}|\vx_{k})),\\[3mm]
\nonumber for~i&=&1,2,\ldots,l_2^\prime+1,~j=1,2,\ldots,l_2^\prime+1.
\end{eqnarray}
$\mD_k^{11}$, $\mD_k^{12}$, $\mD_k^{21}$ and $\mD_k^{22}$ in (\ref{Eqs_sec3_19}) are calculated as follows
\begin{eqnarray}
\label{Eqs_sec3_21} \mD_k^{22}&=&\mB_k^{l_2^\prime+1,l_2^\prime+1}+\mC_k^{1,1},\\
\label{Eqs_sec3_22} \mD_k^{21}&=&\left(
\begin{array}{ccc}
\mB_k^{l_2^\prime+1,1} & \cdots & \mB_k^{l_2^\prime+1,l_2^\prime}\\
\end{array}
\right)=(\mD_k^{12})^T,\\
\label{Eqs_sec3_23} \mD_k^{11}&=&\left(\begin{array}{ccc}
\mB_k^{1,1} &\cdots & \mB_k^{1,l_2^\prime} \\
\vdots & \ddots & \vdots \\
\mB_k^{l_2^\prime,1} & \cdots & \mB_k^{l_2^\prime, l_2^\prime} \\
\end{array}
\right).
\end{eqnarray}
\end{corollary}

Proof: See the Appendix.

\begin{corollary}\label{coro_4}
If the measurement noise is $l$-step auto-correlated ($l\geq0$), the process noise is temporally independent, and the process noise and the
measurement noise are mutually independent, i.e., $l_1=l$, $l_2=0$, $l_3=0$, $l_4=0$, then the sequence $\{\mJ_k\}$ of posterior information
submatrices for estimating state vector $\{\vx_k\}$ obeys the recursion
\begin{eqnarray}
\label{Eqs_sec3_24} \mJ_{k+1}=\mD_k^{22}-\mD_k^{21}(\mD_k^{11}+\mJ_k)^{-1}\mD_k^{12}
\end{eqnarray}
where
\begin{eqnarray}
\nonumber\mD_k^{11}&=&E(-\bigtriangleup_{\vx_k}^{\vx_k}\ln~p(\vx_{k+1}|\vx_k)),\\
\nonumber\mD_k^{21}&=&E(-\bigtriangleup_{\vx_{k+1}}^{\vx_k}\ln~p(\vx_{k+1}|\vx_k))=(\mD_k^{12})^T,\\
\nonumber\mD_k^{22}&=&E(-\bigtriangleup_{\vx_{k+1}}^{\vx_{k+1}}\ln~p(\vx_{k+1}|\vx_k))+\mC_k,
\end{eqnarray}
\begin{eqnarray}
\nonumber &&\mC_{k}=\\
\nonumber&&\left\{
\begin{array}{ll}
E(-\bigtriangleup_{\vx_{k+1}}^{\vx_{k+1}}\ln~p(\vz_{k+1}|\vx_{k+1},\\
\qquad\qquad\qquad\qquad\qquad\vz_k,\ldots,\vz_{k-l+1})), & \hbox{$l\geq1$}\\
E(-\bigtriangleup_{\vx_{k+1}}^{\vx_{k+1}}\ln~p(\vz_{k+1}|\vx_{k+1})), & \hbox{$l=0$}.
\end{array}
\right.
\end{eqnarray}
\end{corollary}

Proof: See the Appendix.

\begin{remark}
When the process noise and the measurement noise are mutually independent and temporally independent, respectively, i.e., $l_1=0$, $l_2=0$,
$l_3=0$, $l_4=0$. Based on Corollary \ref{coro_4}, it is easy to see that the recursion is the same as Proposition 1 in
\cite{Tichavsky-Muravchik-Nehorai98}.
\end{remark}

\section{Numerical Examples}\label{sec_4}
In this section, we consider two  target tracking examples when noises of dynamic systems are temporally correlated.
We compare the new PCRB to already existing
techniques  which include the method of
 \cite{Tichavsky-Muravchik-Nehorai98}, the pre-whitening method, the state augmentation method and the unbiased measurement conversion method given in \cite{BarShalom-Li-Kirubarajan01}. Moreover, based on the PCRB, we can
consider a sensor selection problem, i.e., determine how many sensors should be selected to obtain a desired tracking performance (see, e.g.,
\cite{Joshi-Boyd09,Shen-Varshney14}).

\subsection{Example 1}
Consider a discrete time second order kinematic system driven by temporally correlated noises. This ``correlated noise acceleration model" can
be used in maneuvering tracking \cite{Li-Jilkov03, BarShalom02}.
The discrete time state equation is
\begin{eqnarray}
\label{Eqs_sec4_1}\vx_{k+1}=\left(
\begin{array}{cc}
1 & T \\
0 & 1 \\
\end{array}
\right)\vx_k+\omega_k,
\end{eqnarray}
where the process noise is an one-step correlated moving-average model, i.e.,
\begin{eqnarray}
\label{Eqs_sec4_2}
 \omega_k =\tilde{\omega}_k+ 0.2\tilde{\omega}_{k-1}
\end{eqnarray}
$\{\tilde{\omega}_k\}$ is a Gaussian white noise with zero mean and variance $\tilde{\mQ}_k=\left(
\begin{array}{cc}
\frac{T^3}{3}& \frac{T^2}{2}\\
\frac{T^2}{2}& T \\
\end{array}
\right)$q,  with power spectral density  $q=10$ and sampling interval $T=2$.

The measurement is given by
\begin{eqnarray}
\label{Eqs_sec4_3}\vz_k&=&\left(
\begin{array}{cc}
1 & 0 \\
0& 1\\
\end{array}
\right)\vx_k+\nu_k,
\end{eqnarray}
where measurement noise is considered one-step correlated and one-step cross-correlated with process noise as discussed in
\cite{Mazor-Averbuch-BarShalom-Dayan98,Simon06}, i.e.,
\begin{eqnarray}
\label{Eqs_sec4_4}
 \nu_k=\tilde{\nu}_k+ 0.2\tilde{\nu}_{k-1}+\omega_{k-1},
\end{eqnarray}
where $\{\tilde{\nu}_k\}$ is a Gaussian white noise with zero mean and variance $\tilde{\mR}_k=
\left(\begin{array}{cc}
20^2&0\\
0&5^2\\
\end{array}
\right)
$; $\{\tilde{\nu}_k\}$ and $\{\tilde{\omega}_{k}\}$ are
mutually independent.

By (\ref{Eqs_sec4_1})-(\ref{Eqs_sec4_4}), it can easily be shown that $\{\nu_k\}$ is one-step correlated, $\{\omega_k\}$ is one-step correlated,
$\{\nu_k\}$ is backward two-step and forward one-step correlated with $\{\omega_k\}$ , i.e.,  $l_1=1$, $l_2=1$, $l_3=2$, $l_4=1$.Thus, Theorem \ref{thm_1} can be evaluated.

\begin{figure}[htbp]
\centering \resizebox{12cm}{8cm}{\includegraphics{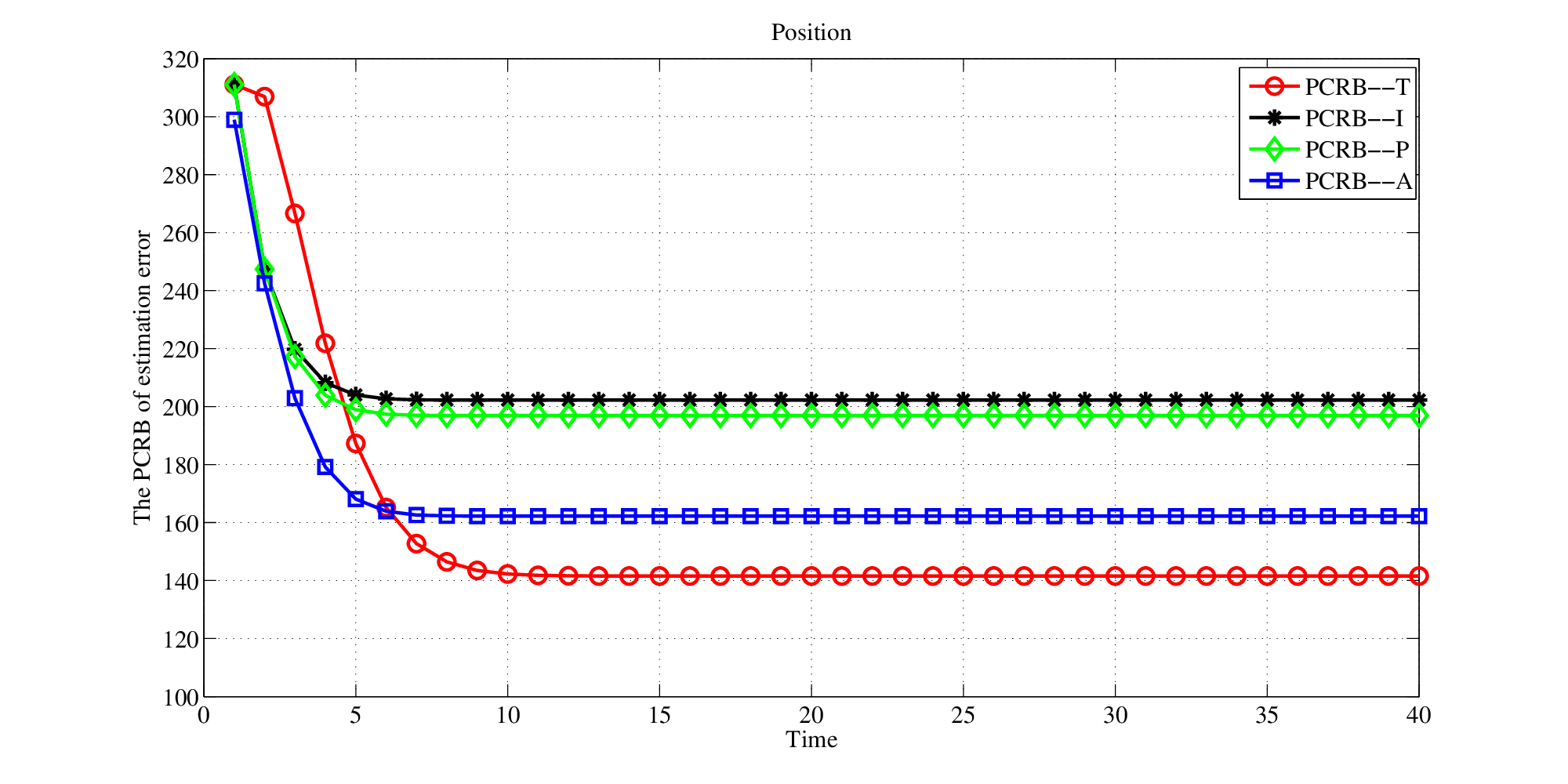}} \caption{The PCRB of the position is plotted as a function of the time step.}\label{fig_1}
\end{figure}

\begin{figure}[htbp]
\centering \resizebox{12cm}{8cm}{\includegraphics{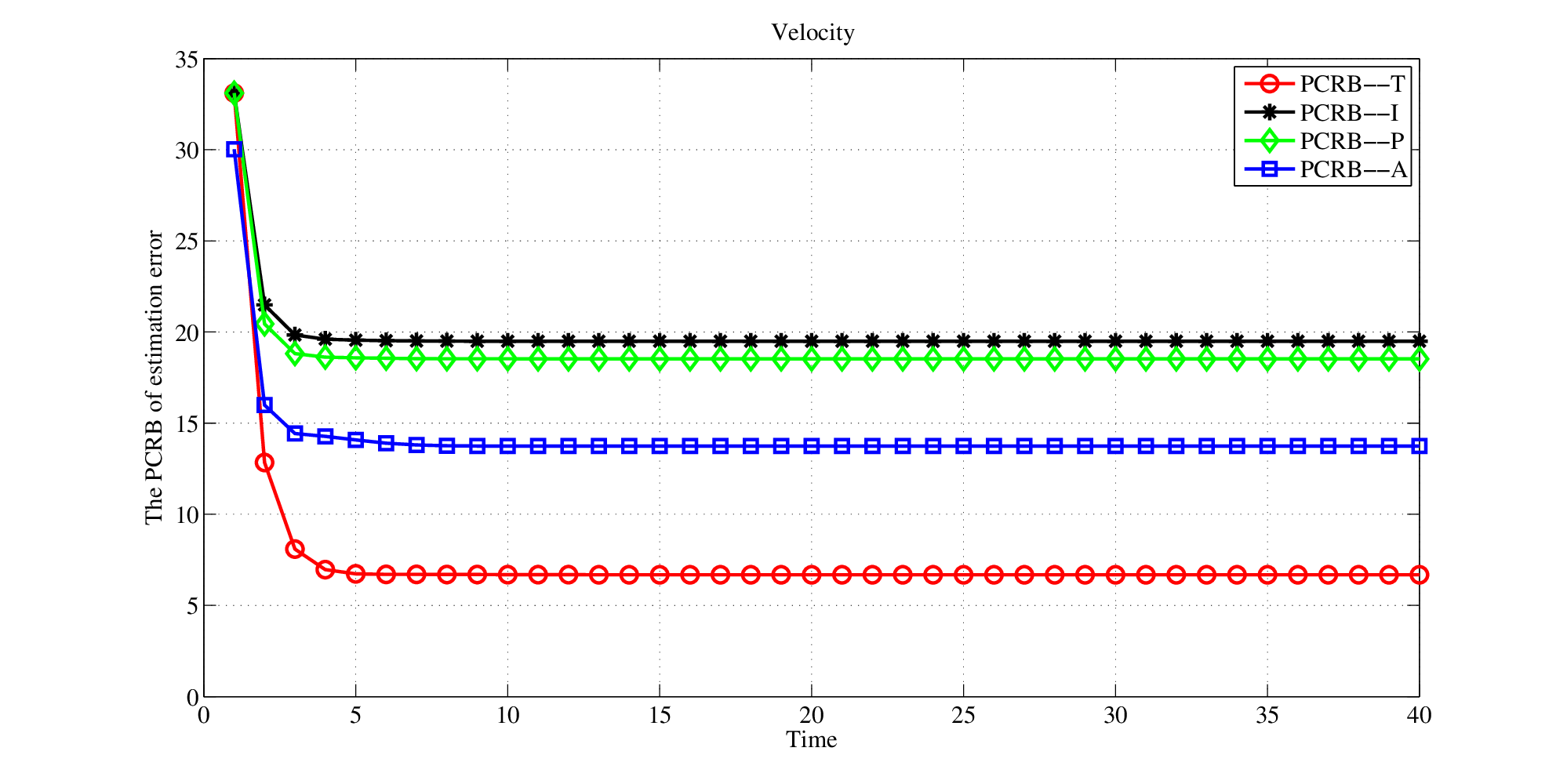}} \caption{The PCRB of the velocity is plotted as a function of the time step.}\label{fig_2}
\end{figure}

\begin{figure}[htbp]
\centering \resizebox{12cm}{8cm}{\includegraphics{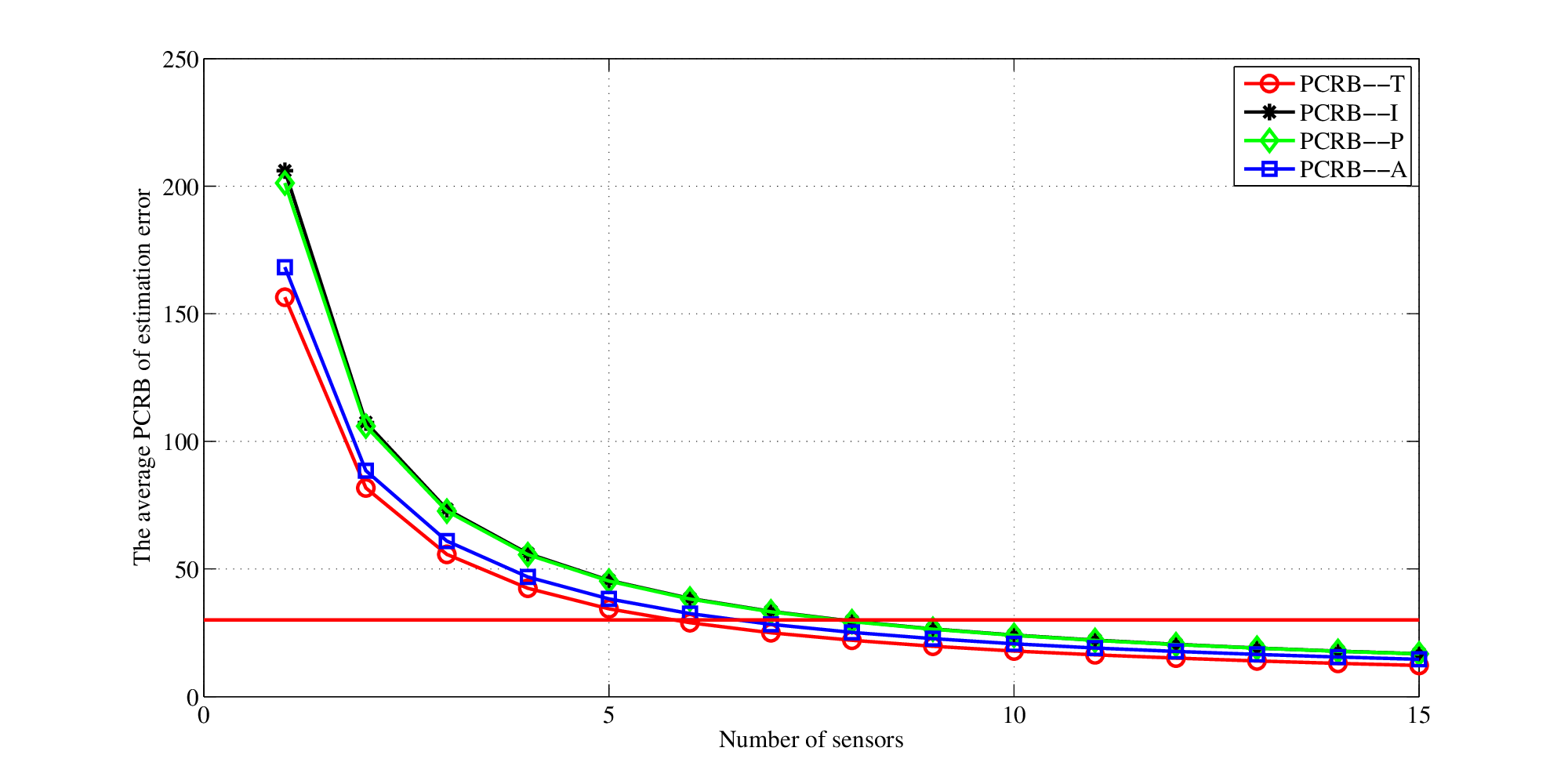}} \caption{The average PCRB of 40 time step of the position is plotted as a function
of number of selected sensors.}\label{fig_3}
\end{figure}

From these assumptions, the conditional probability densities are given as
\begin{eqnarray}
\nonumber &&-\ln~p(\vx_{k+1}|\vx_k,\vz_k)\\
\nonumber &&=c_1+\frac{1}{2}[\vx_{k+1}-g_k(\vx_k,\vz_k)]^T \tilde{\mQ}_k^{-1}[\vx_{k+1}-g_k(\vx_k,\vz_k)],\\[3mm]
\nonumber &&-\ln~p(\vz_{k+1}|\vx_{k+1},\vz_k,\vx_k)\\
\nonumber&&=c_2+\frac{1}{2}[\vz_{k+1}-e_k(\vx_{k+1},\vz_k,\vx_k)]^T\\
\nonumber&&\cdot\tilde{\mR}_{k+1}^{-1}[\vz_{k+1}-e_k(\vx_{k+1},\vz_k,\vx_k)],
\end{eqnarray}
where $c_1$ and $c_2$ are constants, and
\begin{eqnarray}
\nonumber g_k(\vx_k,\vz_k)&=&\left(
\begin{array}{cc}
0.8 & T \\
0 & 0.8 \\
\end{array}
\right)\vx_k+0.2\vz_k\\
\nonumber&&-0.2\tilde{\nu}_k-0.2^2\tilde{\nu}_{k-1}-0.2^2\tilde{\omega}_{k-2}
,\\
\nonumber e_k(\vx_{k+1},\vz_k,\vx_k)&=&\left(
\begin{array}{cc}
2 & 0 \\
0& 2\\
\end{array}
\right)\vx_{k+1}-\left(
\begin{array}{cc}
1.2 &T \\
0& 1.2\\
\end{array}
\right)\vx_k\\
\nonumber&&+0.2\vz_k-0.2^2\tilde{\nu}_{k-1}-0.2\omega_{k-1}.\\
\end{eqnarray}
Using (\ref{Eqs_sec3_3})-(\ref{Eqs_sec3_4}), we can get
\begin{eqnarray}
\nonumber\mB_k^{11}&=&\left(
\begin{array}{cc}
0.8 & T \\
0 & 0.8 \\
\end{array}
\right)^T\tilde{\mQ}_k^{-1}\left(
\begin{array}{cc}
0.8 & T \\
0 & 0.8 \\
\end{array}
\right),\\
\nonumber\mC_k^{11}&=& \left(\begin{array}{cc}
1.2 &T \\
0& 1.2\\
\end{array}
\right)^T\tilde{\mR}_{k+1}^{-1} \left(\begin{array}{cc}
1.2 &T \\
0& 1.2\\
\end{array}
\right),\\
\nonumber\mB_k^{12}&=&- \left(
\begin{array}{cc}
0.8 & T \\
0 & 0.8 \\
\end{array}
\right)^T \tilde{\mQ}_k^{-1},\\
\nonumber\mC_k^{12}&=&-\left(\begin{array}{cc}
1.2 &T \\
0& 1.2\\
\end{array}
\right)^T\tilde{\mR}_{k+1}^{-1} \left(\begin{array}{cc}
2 &0 \\
0&2\\
\end{array}
\right),\\
\nonumber \mB_k^{22}&=&\tilde{\mQ}_k^{-1},~ \mC_k^{22}=\left(\begin{array}{cc}
2 &0 \\
0&2\\
\end{array}
\right)^T\tilde{\mR}_{k+1}^{-1} \left(\begin{array}{cc}
2 &0 \\
0&2\\
\end{array}
\right).
\end{eqnarray}
A straightforward calculation of (\ref{Eqs_sec3_12})-(\ref{Eqs_sec3_14}) gives
\begin{eqnarray}
\nonumber\mD_k^{11}=\mB_k^{11}+\mC_k^{11},\\
\nonumber\mD_k^{12}=\mB_k^{12}+\mC_k^{12},\\
\nonumber\mD_k^{22}=\mB_k^{22}+\mC_k^{22}.
\end{eqnarray}
Thus, we can derive the PCRB for estimating state vector $\{\vx_k\}$ by (\ref{Eqs_sec3_1}) and (\ref{Eqs_sec3_11})-(\ref{Eqs_sec3_14}) of Theorem \ref{thm_1}. The corresponding PCRB is denoted by PCRB-T in Figures \ref{fig_1}--\ref{fig_3}.

Since there are no existing techniques to handle auto-correlation and cross-correlation simultaneously, we approximatively use the state augmentation method, the pre-whitening method, and the method of
 \cite{Tichavsky-Muravchik-Nehorai98}.

For the state augmentation method, we consider auto-correlation of measurement noise
and state noise simultaneously but ignoring the cross-correlation. By (\ref{Eqs_sec4_2}), we can easily get $\omega_k=0.2\omega_{k-1}+\tilde{\omega}_{k}-0.2^2\tilde{\omega}_{k-2}$, which is not
an auto-regressive model. If we let
$\beta_k=\tilde{\omega}_{k}-0.2^2\tilde{\omega}_{k-2}$ which is a Gaussian white noise with zero mean and variance $\tilde{\mQ}_k
+0.2^4\tilde{\mQ}_{k-2}$, and assume that $\{\beta_k\}$ are mutually independent, then the state noise can be approximated by an auto-regressive model $\omega_k=0.2\omega_{k-1}+\beta_k$. Similarly, by (\ref{Eqs_sec4_4}), $\nu_k=0.2\nu_{k-1}+\tilde{\nu_k}-0.2^2\tilde{\nu}_{k-2}+\tilde{\omega}_{k-1}-0.2^2\tilde{\omega}_{k-3}$, let
$\gamma_k=\tilde{\nu_k}-0.2^2\tilde{\nu}_{k-2}+\tilde{\omega}_{k-1}-0.2^2\tilde{\omega}_{k-3}$ which is a Gaussian white noise with zero mean and variance $\tilde{\mQ}_{k-1}+0.2^4\tilde{\mQ}_{k-3}+\tilde{\mR}_{k}+0.2^4\tilde{\mR}_{k-2}$, and assume that $\{\gamma_k\}$ are mutually independent, then the measurement noise can be approximated by an auto-regressive model $\nu_k=0.2\nu_{k-1}+\gamma_k$. Therefore, based on the auto-regressive models $\omega_k=0.2\omega_{k-1}+\beta_k$ and $\nu_k=0.2\nu_{k-1}+\gamma_k$, we can use the the state augmentation method given in \cite{BarShalom-Li-Kirubarajan01} to derive an approximate PCRB which is denoted by PCRB-A in Figures \ref{fig_1}--\ref{fig_3}.

For the pre-whitening method, we consider cross-correlation of measurement noise
and state noise but ignoring auto-correlation of them.  Therefore,  we can use  the pre-whitening method given in \cite{BarShalom-Li-Kirubarajan01} to derive an approximate PCRB which is denoted by PCRB-P in Figures \ref{fig_1}--\ref{fig_3}.

For the method of \cite{Tichavsky-Muravchik-Nehorai98}, we ignore the correlation of noises and assume independent noises. Thus, we can use the method of \cite{Tichavsky-Muravchik-Nehorai98} to derive an approximate PCRB which is denoted by PCRB-I in Figures \ref{fig_1}--\ref{fig_3}.

PCRBs of the position and velocity state are plotted as a function of the time step in Figures \ref{fig_1}--\ref{fig_2}, respectively. For
sensor selection, the average PCRB of 40 time steps is plotted as a function of number of selected sensors in Figure \ref{fig_3}.

The Figures \ref{fig_1}-\ref{fig_2} show that the new PCRB
is  significantly different from those of the other methods. The reason maybe that the approximation loss of the augmentation method and the pre-whitening
method which ignore parts of correlation of noises and cannot deal with
auto-correlation and cross-correlation simultaneously.
In addition, Figures \ref{fig_1}-\ref{fig_2} show that the time-invariant character of the kinematic model implies that the PCRB converges to a constant after some time steps.
In Figure \ref{fig_3}, it can be seen that when the number of selected sensors is increasing, the gap of PCRB becomes smaller. Figure
\ref{fig_3} also shows that if we want to achieve PCRB of the estimation error less than 30 $m$, 6 sensors have to be used based on the new PCRB at least. However, if we only consider the case of the auto-correlation of the state and measurement noises, 7 sensors have to be used, the other cases may be used more than 8 sensors. Thus, number of selected sensors becomes very different to obtain a desired estimation performance.


\subsection{Example 2}

In this example, we consider a discrete time dynamic system  with \emph{nonlinear} measurements as follows. The four-dimensional state variable
includes position and velocity $(x, \dot{x}, y, \dot{y})$ driven by correlated noise, respectively,
\begin{eqnarray}%
\label{Eqs_sec4_5}&&\vx_{k+1}= \mF_k\vx_k+\omega_k, ~~~ \mF_k=\left(
                                 \begin{array}{cccc}
                                   1 & T & 0 & 0 \\
                                   0 & 1& 0 &0 \\
                                   0 & 0 & 1 & T \\
                                   0 & 0 & 0& 1 \\
                                 \end{array}
                               \right)
\end{eqnarray}
where the process noise is an two-step correlated moving-average model,
\begin{eqnarray}
\label{Eqs_sec4_6} \omega_k =\tilde{\omega}_k+ \tilde{\omega}_{k-1}+ \tilde{\omega}_{k-2}.
\end{eqnarray}
$\{\tilde{\omega}_k\}$ is a Gaussian white noise with zero mean and variance
\begin{eqnarray}%
\label{Eqs_sec4_7}&&\tilde{\Sigma}_1= \left(
                                 \begin{array}{cccc}
                                   \frac{T^3}{3} & \frac{T^2}{2} & 0 & 0 \\
                                   \frac{T^2}{2} & T&0 &0 \\
                                   0 & 0 &\frac{T^3}{3} & \frac{T^2}{2}\\
                                   0 & 0 & \frac{T^2}{2}& T \\
                                 \end{array}
                               \right)q,
\end{eqnarray}
with sampling interval $T=3$ and power spectral density $q=10$.

The two-dimensional nonlinear measurement vector includes range and azimuth, respectively,
\begin{eqnarray}
\label{Eqs_sec4_8}\vz_k=\vh_k(\vx_k)+\vv_k,
\end{eqnarray}
where the nonlinear measurement function is
\begin{eqnarray}
\label{Eqs_sec4_9}\vh_k(\vx_k)=\left(
\begin{array}{c}
\vh_k^1(\vx_k) \\
\vh_k^2(\vx_k) \\
\end{array}
\right)=\left(
\begin{array}{c}
\sqrt{(\vx_k^1)^2+(\vx_k^3)^2} \\
\tan^{-1}(\frac{\vx_k^3}{\vx_k^1}) \\
\end{array}
\right).
\end{eqnarray}
$\vx_k^i$ is the $i$-th entry of the state vector $\vx_k$. $\{\vv_k\}$ is a Gaussian white noise with zero mean and variance matrix
\begin{eqnarray}%
\label{Eqs_sec4_10}&&\Sigma_2= \left(
                                 \begin{array}{cc}
                                   50^2 & 0  \\
                                   0 &  0.01^2\\
                                  \end{array}
                               \right).
\end{eqnarray}
$\{\vv_k\}$ and $\{\tilde{\omega}_k\}$ are mutually independent.

\begin{figure}[htbp]
\centering \resizebox{12cm}{8cm}{\includegraphics{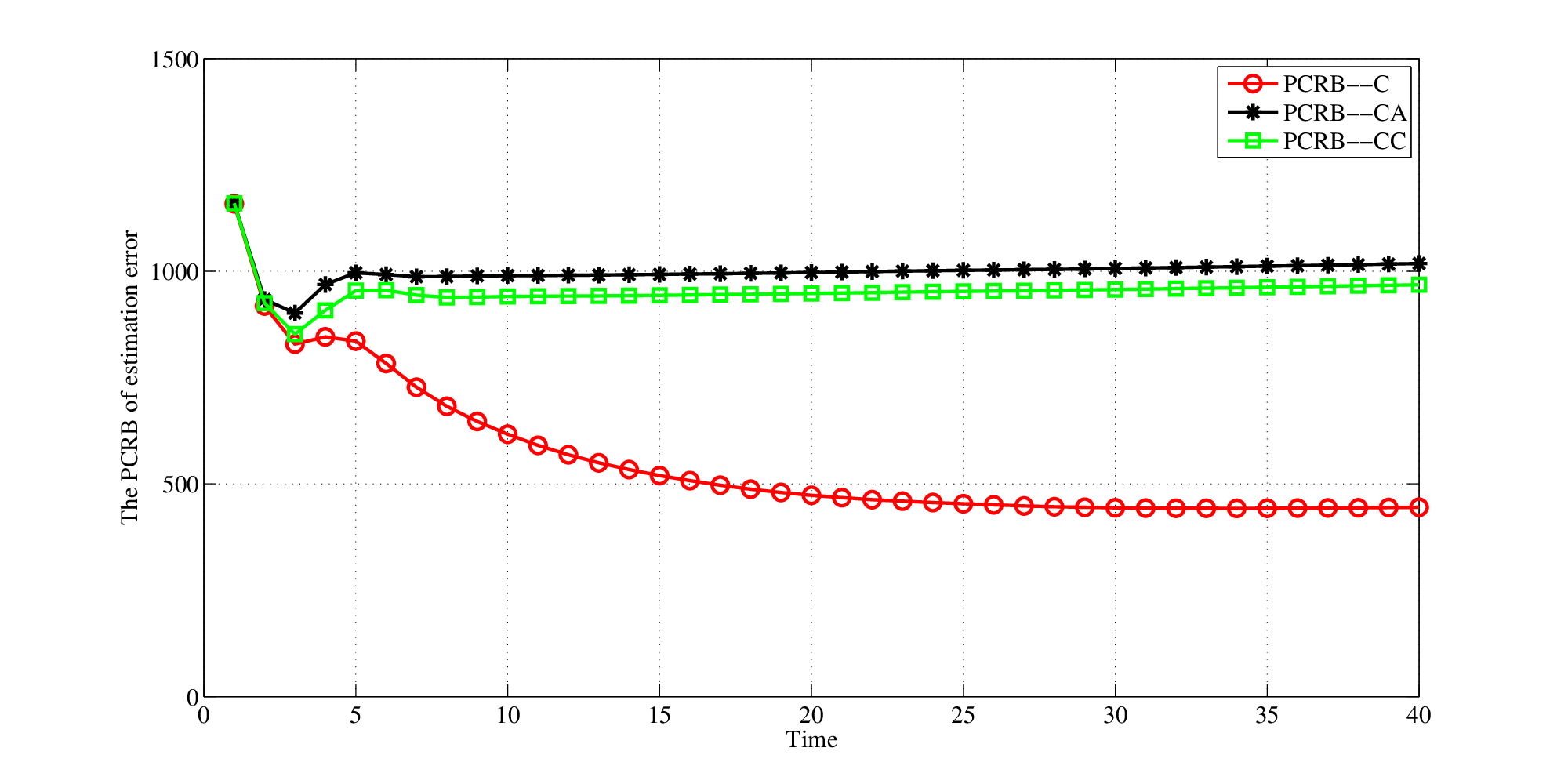}} \caption{The PCRB of the position is plotted as a function of the time step.}\label{fig_4}
\end{figure}

\begin{figure}[htbp]
\centering \resizebox{12cm}{8cm}{\includegraphics{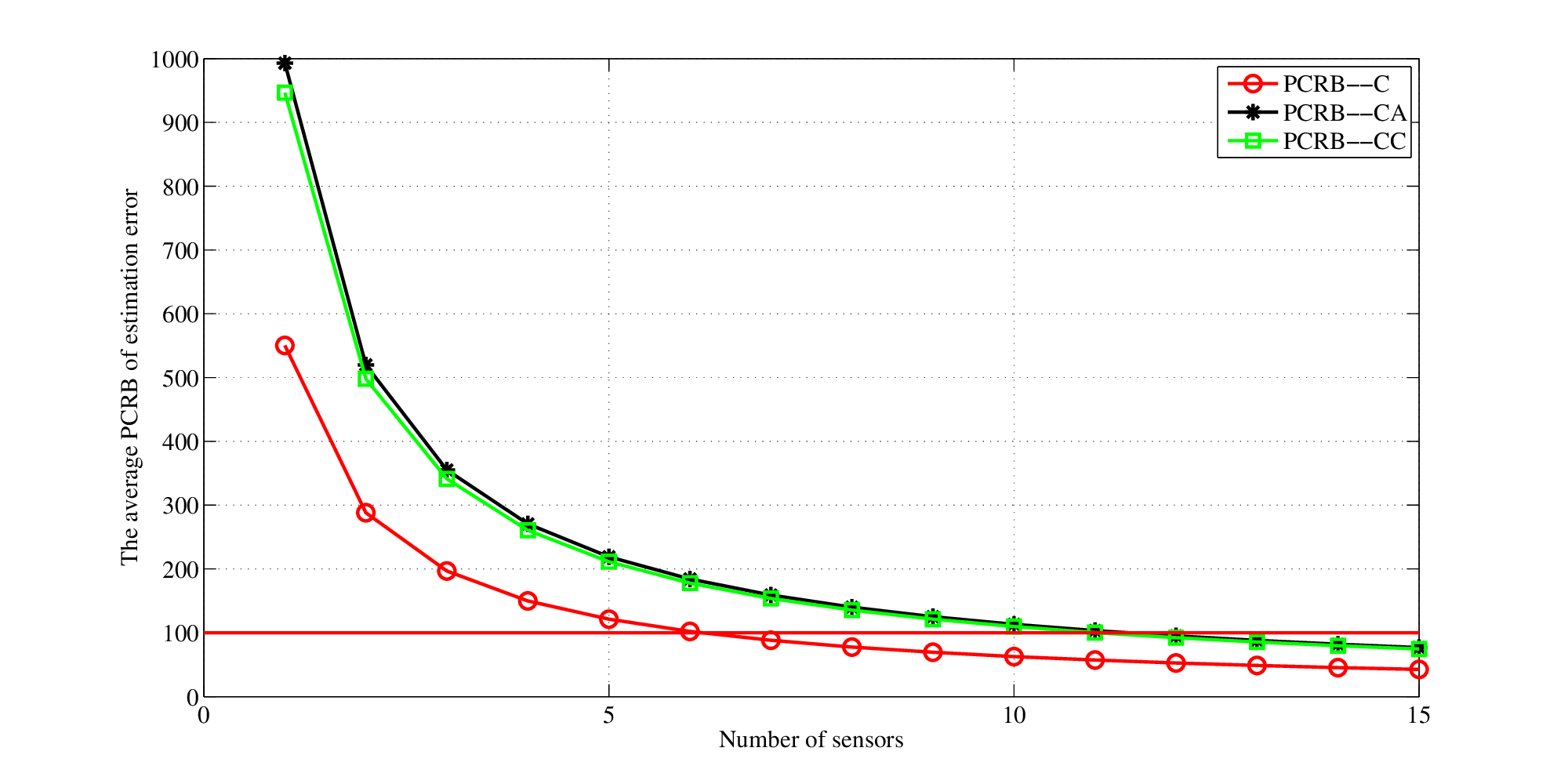}} \caption{The average PCRB of 40 time steps is plotted as a
function of number of selected sensors.}\label{fig_5}
\end{figure}

It can easily be seen that the model satisfies Corollary \ref{coro_3} and $l=2$. From these assumptions, the conditional probability densities
are given as
\begin{eqnarray}
\nonumber &&-\ln p(\vx_{k+1}|\vx_k, \vx_{k-1})\\
\nonumber&&=c_3+(\vx_{k+1}-g(\vx_k, \vx_{k-1}))^T\tilde{\Sigma}_1^{-1}(\vx_{k+1}-g(\vx_k, \vx_{k-1})),\\[3mm]
\nonumber &&-\ln p(\vz_{k+1}|\vx_{k+1})\\
\nonumber&&=c_4+(\vz_{k+1}-h(x_{k+1}))^T\Sigma_2^{-1}(\vz_{k+1}-h(x_{k+1})),\\[3mm]
\nonumber &&g(\vx_k, \vx_{k-1})\\
\nonumber&&=(\mI+\mF_k)\vx_k-\mF_k\vx_{k-1}-\tilde{\omega}_{k-3},
\end{eqnarray}
where $c_3$ and $c_4$ are constants; $\mI$ is an identity matrix with compatible dimensions. A straightforward calculation of
(\ref{Eqs_sec3_019})-(\ref{Eqs_sec3_020}) gives
\begin{eqnarray}
\nonumber \mB_k^{11}&=&\mF_k^T\tilde{\Sigma}_1^{-1}\mF_k,~~ \mB_k^{12}=-\mF_k^T\tilde{\Sigma}_1^{-1}(\mI+\mF_k),\\
\nonumber\mB_k^{13}&=&\mF_k^T\tilde{\Sigma}_1^{-1},~~\mB_k^{21}=(\mB_k^{12})^T,\\
\nonumber\mB_k^{22}&=&(\mI+\mF_k)^T\tilde{\Sigma}_1^{-1}(\mI+\mF_k),~~ \mB_k^{23}=-(\mI+\mF_k)^T\tilde{\Sigma}_1^{-1},\\
\nonumber \mB_k^{31}&=&(\mB_k^{13})^T,~~ \mB_k^{32}=(\mB_k^{23})^T, ~~\mB_k^{33}=\tilde{\Sigma}_1^{-1},\\
\nonumber \mC_k^{11}&=&E\{[\nabla_{\vx_{k+1}}{h(\vx_{k+1})}^T]\Sigma_2^{-1}[\nabla_{\vx_{k+1}}{h(\vx_{k+1})}^T]^T\}.
\end{eqnarray}
$\mC_k^{11}$ can be calculated by numerical Monte-Carlo methods. Using (\ref{Eqs_sec3_21})-(\ref{Eqs_sec3_23}), we can easily get
\begin{eqnarray}
\label{Eqs_sec4_07}\mD_k^{11}&=&\left(
\begin{array}{cc}
\mB_k^{11} & \mB_k^{12} \\
\mB_k^{21}& \mB_k^{22}\\
\end{array}
\right),\\
\nonumber \mD_k^{21}&=&\left(
\begin{array}{cc}
\mB_k^{31} & \mB_k^{32} \\
\end{array}
\right)=(\mD_k^{12})^T,\\
\nonumber \mD_k^{22}&=&\mB_k^{33}+\mC_k^{11}.
\end{eqnarray}
Combing (\ref{Eqs_sec4_07}), (\ref{Eqs_sec3_19}) and (\ref{Eqs_sec3_20}), after some simplification, we have the simpler recursion
\begin{eqnarray}
\nonumber \mJ_{k+1}&=&\mD_k^{22}-\mD_k^{21}(\mE_k+\mD_k^{11})^{-1}\mD_k^{12},\\
\nonumber \mE_k^{11}&=&\mE_{k-1}^{22}+\mB_{k-1}^{22}-(\mE_{k-1}^{21}+\mB_{k-1}^{21})\\
\nonumber&&\cdot(\mE_{k-1}^{11}+\mB_{k-1}^{11})^{-1}(\mE_{k-1}^{12}+\mB_{k-1}^{12}),\\
\nonumber
\mE_k^{12}&=&\mB_{k-1}^{23}-(\mE_{k-1}^{21}+\mB_{k-1}^{21})(\mE_{k-1}^{11}+\mB_{k-1}^{11})^{-1}\mB_{k-1}^{13}\\
\nonumber&=&(\mE_k^{21})^T,\\
\nonumber \mE_k^{22}&=&\mB_{k-1}^{33}+\mC_{k-1}^{11}-\mB_{k-1}^{31}(\mE_{k-1}^{11}+\mB_{k-1}^{11})^{-1}\mB_{k-1}^{13}.
\end{eqnarray}
Thus, based on the above recursion, we can derive the PCRB for estimating state vector $\{\vx_k\}$ by Corollary \ref{coro_3}. The corresponding PCRB is denoted by PCRB-C in Figures \ref{fig_4}--\ref{fig_5}. 

Since there are no existing methods to deal with this example accurately, we approximatively use the unbiased measurement conversion method given in
 \cite{BarShalom-Li-Kirubarajan01}, which can convert the nonlinear system into linear system. Moreover, similar to Example 1, we use the state augmentation method to derive an approximate PCRB which is denoted by PCRB-CA in Figures \ref{fig_4}--\ref{fig_5}. In addition, for the converted linear system, we can also use Corollary \ref{coro_3} to derive an approximate PCRB which is denoted by PCRB-CC in Figures \ref{fig_4}--\ref{fig_5}. 


In Figure \ref{fig_4}, PCRB of the position  is plotted as a function of the time step. For
sensor selection, the \emph{average} PCRB of 40 time steps is plotted as a function of number of selected sensors in Figure \ref{fig_5}.

Figures \ref{fig_4}-\ref{fig_5} shows that  the new PCRB
is  significantly different from the other two methods. The reason maybe the approximation loss of unbiased conversion, and that the
pdf of the noise of the converted linear system is non-Gaussian and uncertain, and the approximation loss of the augmentation method which ignore parts of correlation of noises.
Figure \ref{fig_5} also shows that if we want to achieve PCRB of the estimation error less than 100 $m$, 6 sensors have to used at least. However, if we use the other methods, we have to select 11 sensors.  Thus, this example also shows that number of selected sensors becomes very different to obtain a desired estimation performance.


\section{Conclusions}\label{sec_5}
In this paper, we have derived  a unified recursive formula of the mean-square error lower bound for the discrete-time nonlinear filtering
problem when noises of dynamic systems are temporally correlated based on the posterior version of the Cram$\acute{e}$r-Rao
inequality. It can be applicable to the multi-step correlated process noise, multi-step correlated measurement noise and multi-step
cross-correlated process and measurement noise simultaneously. Although the unified formula can come across the three cases of finite-step
correlated noises, a few corollaries with simpler formulae follow to elucidate special cases, which may be used more frequently. Two
typical target tracking examples have shown that the new PCRB is significantly different from that of the other existing approximation methods which include
the method of ignoring correlation of noises, the pre-whitening method, the state augmentation method and the unbiased measurement conversion method. Thus, when they are applied
to sensor selection problems, simulations show that number of selected sensors becomes very different to obtain a desired estimation performance. Future research challenges include sensor management when noises of dynamic systems are temporally correlated for
multitarget tracking and data association.

\section{Appendix}\label{sec_6}
\begin{lemma}\label{lemma_1}
(see, e.g., \cite{Horn-Johnson12,Jiang-Zhou-Zhu10}) Consider a partitioned matrix
\begin{eqnarray}
\nonumber\mA=\left(
\begin{array}{cc}
\mA_{11} & \mA_{12} \\
\mA_{21} & \mA_{22} \\
\end{array}
\right).
\end{eqnarray}
If $\mA$ is invertible, then $\mA_{11}$ is invertible if and only if the Schur complement $\triangle=\mA_{22}-\mA_{21}\mA_{11}^{-1}\mA_{12}$ is
invertible, and
\begin{eqnarray}
\nonumber\mA^{-1}&=&\left(
\begin{array}{cc}
I & -\mA_{11}^{-1}\mA_{12} \\
0 & I \\
\end{array}
\right)\left(
\begin{array}{cc}
\mA_{11}^{-1} & 0 \\
0 & \triangle^{-1}\\
\end{array}
\right)\\
\nonumber&&\cdot\left(
\begin{array}{cc}
I & 0 \\
-\mA_{21}\mA_{11}^{-1} & I \\
\end{array}
\right).
\end{eqnarray}
\end{lemma}

\begin{lemma}\label{lemma_2}
(see, e.g., \cite{Horn-Johnson12,Jiang-Zhou-Zhu10}) Let $\mB=\left(
\begin{array}{cc}
\mB_{11} & \mB_{12} \\
\end{array}
\right)$, ~$\mA=\left(
\begin{array}{cc}
\mA_{11} & \mA_{12} \\
\mA_{21} & \mA_{22} \\
\end{array}
\right)$,~ $\mC=\left(
\begin{array}{c}
\mC_{11} \\
\mC_{21} \\
\end{array}
\right)$. If $\mA$ and $\mA_{11}$ are both invertible, then $\triangle=\mA_{22}-\mA_{21}\mA_{11}^{-1}\mA_{12}$ is invertible, and use the Lemma
\ref{lemma_1}, we have
\begin{eqnarray}
\nonumber\mB\mA^{-1}\mC&=&\mB_{11}\mA_{11}^{-1}\mC_{11}+(\mB_{12}-\mB_{11}\mA_{11}^{-1}\mA_{12})\triangle^{-1}\\
&&\cdot(\mC_{21}-\mA_{21}\mA_{11}^{-1}\mC_{11}).
\end{eqnarray}
\end{lemma}

\subsection{proof of Proposition \ref{thm_1}}
According to the three definitions of \emph{finite-step correlated} noises given in Section \ref{sec_2_1} and
(\ref{Eqs_sec2_1})-(\ref{Eqs_sec2_2}), to derive a recursive formula,  $p(\vx_{k+1}|\mX_k,\mZ_k)$ and $p(\vz_{k+1}|\mX_{k+1},\mZ_k)$ of
(\ref{Eqs_sec2_8}) may be approximately written as
\begin{flushleft}
$p(\vx_{k+1}|\mX_k,\mZ_k)\approx$
\end{flushleft}
\begin{eqnarray}
\label{Eqs_secA_1} \left\{
\begin{array}{ll}
p(\vx_{k+1}|\vx_k) & \hbox{for ~$l_2=0$ or $l_2=1$} \\
~&\hbox{ and~$l_4=0$ ~}\\
p(\vx_{k+1}|\vx_k,\ldots,\vx_{k-l_2+1})& \hbox{for ~ $l_2\geq2$,~$l_4=0$}\\
p(\vx_{k+1}|\vx_k,\\
\qquad~~~\vz_k,\ldots,\vz_{k-l_4+1})& \hbox{ for $l_2=0$}\\
~&\hbox{or $l_2=1$~and $l_4>0$}\\
p(\vx_{k+1}|\vx_k,\ldots,\vx_{k-l_2+1},\\
\qquad~~~\vz_k,\ldots,\vz_{k-l_4+1})& \hbox{ for ~$l_2\geq2$,~$l_4>0$},\\
\end{array}
\right.
\end{eqnarray}
\begin{flushleft}
$p(\vz_{k+1}|\mX_{k+1},\mZ_k)\approx$
\end{flushleft}
\begin{eqnarray}
\label{Eqsm_31}\left\{
\begin{array}{ll}
p(\vz_{k+1}|\vx_{k+1}) & \hbox{for ~ $l_3=0$ or $l_3=1$} \\
~&\hbox{and~$l_1=0$}\\
p(\vz_{k+1}|\vx_{k+1},\ldots,&\hbox{for ~ $l_3\geq2$,}\\
\qquad\qquad\vx_{k-l_3+2})& \hbox{$l_1=0$}\\
p(\vz_{k+1}|\vx_{k+1},& \hbox{ for ~ $l_3=0$}\\
\qquad~~~ \vz_k,\ldots,\vz_{k-l_1+1})&\hbox{or $l_3=1$ and~$l_1>0$}\\
p(\vz_{k+1}|\vx_{k+1},\ldots,\vx_{k-l_3+2}\\
\qquad~~~\vz_k,\ldots,\vz_{k-l_1+1})& \hbox{ for ~$l_3\geq2$,~$l_1>0$}.\\
\end{array}
\right.
\end{eqnarray}
If we denote by $l_3^\prime\triangleq \max\{l_3,1\},l_2^\prime\triangleq max\{l_2,1\}$, then (\ref{Eqs_secA_1}) and (\ref{Eqsm_31}) can be
simplified as
\begin{eqnarray}\label{Eqsm_32}
\nonumber&&p(\vx_{k+1}|\mX_k,\mZ_k)\\
&&\approx\left\{
\begin{array}{ll}
p(\vx_{k+1}|\vx_k,\ldots,\vx_{k-l_2^\prime+1})& \hbox{for ~ $l_4=0$}\\
p(\vx_{k+1}|\vx_k,\ldots,\vx_{k-l_2^\prime+1},\\
\qquad~~~\vz_k,\ldots,\vz_{k-l_4+1})& \hbox{ for ~$l_4>0$},\\
\end{array}
\right.
\end{eqnarray}
\begin{eqnarray}\label{Eqsm_33}
\nonumber&&p(\vz_{k+1}|\mX_{k+1},\mZ_k)\\
&&\approx \left\{
\begin{array}{ll}
p(\vz_{k+1}|\vx_{k+1},\ldots,\vx_{k-l_3^\prime+2})& \hbox{for ~ $l_1=0$}\\
p(\vz_{k+1}|\vx_{k+1},\ldots,\vx_{k-l_3^\prime+2},\\
\qquad~~~\vz_k,\ldots,\vz_{k-l_1+1})& \hbox{ for ~$l_1>0$}.\\
\end{array}
\right.
\end{eqnarray}

In order to derive the recursion of $\mJ_{k+1}$, we needs to decompose vector $\mX_k$ and $\mX_{k+1}$ . Based on
(\ref{Eqsm_32})-(\ref{Eqsm_33}), we find that the decomposition depends on $l_3^\prime$ and $l_2^\prime+1$. Thus, we discuss three cases
$l_3^\prime>l_2^\prime+1$, $l_3^\prime<l_2^\prime+1$ and $l_3^\prime=l_2^\prime+1$, respectively.
%
\begin{enumerate}
\item $l_3^\prime>l_2^\prime+1$

We decompose $\mX_k$ as
\begin{eqnarray}
 \mX_k=(\mX_{k-l_3^\prime+1}^\prime,\ldots,\vx_{k-l_2^\prime+1}^\prime,\ldots ,\vx_{k-1}^\prime,\vx_k^\prime)^\prime
 \end{eqnarray}
 and $J(\mX_k)$ correspondingly as
 \begin{flushleft}
 $J(\mX_k)=$
 \end{flushleft}
\begin{eqnarray}
\nonumber\left(
\begin{array}{ccccc}
\mA_k^{1,1} & \cdots &\mA_k^{1,l_3^\prime-l_2^\prime}&\cdots& \mA_k^{1,l_3^\prime}\\
\vdots &\ldots& \vdots &\cdots& \vdots \\
\mA_k^{l_3^\prime-l_2^\prime,1}&\cdots&\mA_k^{l_3^\prime-l_2^\prime,l_3^\prime-l_2^\prime}&\cdots&\mA_k^{l_3^\prime-l_2^\prime,l_3^\prime}\\
\vdots &\ldots& \vdots &\cdots& \vdots \\
\mA_k^{l_3^\prime,1} & \cdots & \mA_k^{l_3^\prime,l_3^\prime-l_2^\prime}&\cdots&\mA_k^{l_3^\prime,l_3^\prime}\\
\end{array}
\right)
\end{eqnarray}
where
\begin{eqnarray}\label{Eqsm_40}
\nonumber\mA_k^{1,1}&=&E(-\bigtriangleup_{\mX_{k-l_3^\prime+1}}^{\mX_{k-l_3^\prime+1}}\ln~p_k),\\
\nonumber\mA_k^{1,j}&=&E(-\bigtriangleup_{\mX_{k-l_3^\prime+1}}^{\vx_{k-l_3^\prime+j}}\ln~p_k),\\
\nonumber\mA_k^{i,1}&=&E(-\bigtriangleup_{\vx_{k-l_3^\prime+i}}^{\mX_{k-l_3^\prime+1}}\ln~p_k),\\
\nonumber\mA_k^{i,j}&=&E(-\bigtriangleup_{\vx_{k-l_3^\prime+i}}^{\vx_{k-l_3^\prime+j}}\ln~p_k),\\
  &&for ~i,j=2,\ldots,l_3^\prime.
\end{eqnarray}
Using (\ref{Eqs_sec2_8}), (\ref{Eqsm_32}), (\ref{Eqsm_33}) the posterior information matrix for $\mX_{k+1}$ can be written in block form as
\begin{flushleft}
$J(\mX_{k+1})=$
\end{flushleft}
\begin{eqnarray}
\nonumber\left(
           \begin{array}{cccc}
             \mA_k^{1,1} & \cdots  &\mA_k^{1,l_3^\prime} & 0  \\
             \vdots & \ldots & \ldots & \vdots  \\
             \mA_k^{l_3^\prime,1} & \ldots & \mA_k^{l_3^\prime,l_3^\prime} & \mC_k^{l_{31}^\prime,l_3^\prime} \\
             ~ & ~ & +\mC_k^{l_{31}^\prime,l_{31}^\prime} & +\mB_k^{l_{2}^\prime,l_{21}^\prime} \\
             ~ & ~ &+\mB_k^{l_{2}^\prime,l_{2}^\prime} & ~  \\[3mm]
             0&\ldots&\mC_k^{l_3^\prime,l_{31}^\prime}& \mC_k^{l_3^\prime,l_3^\prime} \\
             ~&~&+\mB_k^{l_{21}^\prime,l_{2}^\prime}&+\mB_k^{l_{21}^\prime,l_{21}^\prime}\\
           \end{array}
         \right)
\end{eqnarray}
where $0$'s stand for the zero blocks of appropriate dimensions; $\mB_k^{i,j}$ and $\mC_k^{i,j}$ are defined in
(\ref{Eqs_sec3_3})-(\ref{Eqs_sec3_4}). To simplify the matrix $J(\mX_{k+1})$, we denote by $l_{32}^\prime\triangleq l_3^\prime-l_2^\prime$,
$l_{321}^\prime\triangleq l_3^\prime-l_2^\prime+1$, $l_{31}^\prime\triangleq l_3^\prime-1$, $l_{21}^\prime\triangleq l_2^\prime+1$.

Moreover, the information submatrix $\mJ_{k+1}$ for estimating $\vx_{k+1}\in R^r$ is given as the inverse of the ($r\times r$) right-lower block
of $[J(\mX_{k+1})]^{-1}$, i.e.,
\begin{eqnarray}
\nonumber&&\mJ_{k+1}= \mC_k^{l_3^\prime,l_3^\prime}+\mB_k^{l_{21}^\prime,l_{21}^\prime}\\
 \nonumber&&- \left(
\begin{array}{cccc}
0&\mC_k^{l_3^\prime,1} & \cdots &\mC_k^{l_3^\prime,l_3^\prime-1}\\
~&~&~&+\mB_k^{l_2^\prime+1,l_2^\prime}\\
\end{array}
\right)\\[3mm]
\nonumber&&\cdot\left(
  \begin{array}{cccc}
    \mA_k^{1,1} & \mA_k^{1,2} & \ldots & \mA_k^{1,l_3^\prime} \\
    \mA_k^{2,1} & \mA_k^{2,2} & \ldots & \mA_k^{2,l_3^\prime} \\
    ~ & +\mC_k^{1,1} & \ldots & +\mC_k^{1,l_{31}^\prime} \\
    \ldots& \ldots & \ldots & \ldots \\
    \mA_k^{l_3^\prime,1} &\mA_k^{l_3^\prime,2} & \ldots & \mA_k^{l_3^\prime,l_3^\prime}\\
~&+\mC_k^{l_{31}^\prime,1}&~&+\mC_k^{l_{31}^\prime,l_{31}^\prime}\\
~ &~ & ~ &+\mB_k^{l_{2}^\prime,l_{2}^\prime} \\
  \end{array}
\right)^{-1}\\[3mm]
&&\cdot \left(
\begin{array}{cccc}
0&\mC_k^{l_3^\prime,1} & \cdots&\mC_k^{l_3^\prime,l_3^\prime-1}\\
~&~&~&+\mB_k^{l_2^\prime+1,l_2^\prime}\\
\end{array}
\right)^T.
\end{eqnarray}

Using Lemma \ref{lemma_2} and (\ref{Eqs_sec3_5})-(\ref{Eqs_sec3_7}), it follows that
\begin{eqnarray}
\label{Eqs_secA_06}\mJ_{k+1}&=& \mD_k^{22}-\mD_k^{21}(\mE_k+\mD_k^{11})^{-1}\mD_k^{12},\\
\label{Eqsm_41}\mE_k^{i,j}&=&\mA_k^{i+1,j+1}-\mA_k^{i+1,1}(\mA_k^{1,1})^{-1}\mA_k^{1,j+1}\\
\nonumber&=&(\mE_k^{j,i})^\prime,\\
\nonumber &&for~i,j=1,\ldots,l_3^\prime-1.
\end{eqnarray}
Since  $\mD_k^{11}$, $\mD_k^{12}$, $\mD_k^{21}$ and $\mD_k^{22}$ in (\ref{Eqs_secA_06}) can be calculated as
(\ref{Eqs_sec3_5})-(\ref{Eqs_sec3_7}), we only need to derive the recursion of $\mE_k^{i,j}$.


Based on the  recursion between $p_k$ and $p_{k-1}$ given in (\ref{Eqs_sec2_8}) which depends on (\ref{Eqsm_32}) and (\ref{Eqsm_33}), and using
the definitions of $\mA_k^{i,j}$, $\mB_k^{i,j}$, $\mC_k^{i,j}$ given in (\ref{Eqsm_40}), (\ref{Eqs_sec3_3}), (\ref{Eqs_sec3_4}), respectively, it
follows that
\begin{flushleft}
$\left(
\begin{array}{cc}
\mA_k^{1,1} & \mA_k^{1,j+1} \\
\mA_k^{i+1,1}& \mA_k^{i+1,j+1} \\
\end{array}
\right)$
\end{flushleft}
\begin{eqnarray}
\nonumber=\left(
\begin{array}{ccc}
\mA_{k-1}^{1,1} & \mA_{k-1}^{1,2} & \mA_{k-1}^{1,j+2} \\[3mm]
\mA_{k-1}^{2,1} & \mA_{k-1}^{2,2}&\mA_{k-1}^{2,j+2} \\
~&+\mC_{k-1}^{11}&+\mC_{k-1}^{1,j+1}\\[3mm]
\mA_{k-1}^{i+2,1} & \mA_{k-1}^{i+2,2} &
\mA_{k-1}^{i+2,j+2}\\
~&+\mC_{k-1}^{i+1,1}&+\mC_{k-1}^{i+1,j+1}\\
~&~&+\mB_{k-1}^{i+2-l_3^\prime+l_2^\prime,j+2-l_3^\prime+l_2^\prime}
\end{array}
\right)
\end{eqnarray}
\begin{eqnarray}
\label{Eqsm_16}
\end{eqnarray}
\begin{center}
$for~i=1, \ldots, l_3^\prime-1, j=1, \ldots, l_3^\prime-1$ .
\end{center}
where if $i\leq l_3^\prime-l_2^\prime-2$ or $j\leq l_3^\prime-l_2^\prime-2$, then
$\mB_{k-1}^{i+l_2^\prime-l_3^\prime+2,j+l_2^\prime-l_3^\prime+2}=0$; if $i>l_3^\prime$ or $j>l_3^\prime$, then $\mA_{k-1}^{i,j}=0$.

Note that the matrix $\mE_k^{ij}$ in (\ref{Eqsm_41}) is the Schur complement  of the block $\mA_k^{1,1}$ of the left matrix of Equation
(\ref{Eqsm_16}) and the Schur complement of the corresponding block of the right matrix of Equation (\ref{Eqsm_16}) is
\begin{eqnarray}
\nonumber&&\mA_{k-1}^{i+2,j+2}+\mC_{k-1}^{i+1,j+1}+\mB_{k-1}^{i+2-l_3^\prime+l_2^\prime,j+2-l_3^\prime+l_2^\prime}\\[3mm]
\label{Eqs_A_61}&-&\left(
  \begin{array}{cc}
    \mA_{k-1}^{i+2,1} & \mA_{k-1}^{i+2,2}+\mC_{k-1}^{i+1,1}  \\
  \end{array}
\right)\\[3mm]
\nonumber&\cdot&\left(
         \begin{array}{cc}
           \mA_{k-1}^{1,1} & \mA_{k-1}^{1,2} \\[3mm]
           \mA_{k-1}^{2,1} & \mA_{k-1}^{2,2}+ \\
           ~&\mC_{k-1}^{11}\\
         \end{array}
       \right)^{-1}\left(
                \begin{array}{c}
                  \mA_{k-1}^{1,j+2} \\[3mm]
                  \mA_{k-1}^{2,j+2}+ \\
               \mC_{k-1}^{1,j+1}\\
                \end{array}
              \right).
\end{eqnarray}
Using Lemma \ref{lemma_2}, we can simplify (\ref{Eqs_A_61}). Moreover, by the definition of $\mE_{k-1}^{ij}$ in (\ref{Eqsm_41}), we have the
recursion of $\mE_k^{i,j}$ as
\begin{eqnarray}
\nonumber
\mE_k^{ij}&=&\mE_{k-1}^{i+1,j+1}+\mB_{k-1}^{i+2-l_3^\prime+l_2^\prime,j+2-l_3^\prime+l_2^\prime}+\mC_{k-1}^{i+1,j+1}\\[3mm]
\label{Eqs_secA_01} &&-(\mE_{k-1}^{i+1,1}+\mC_{k-1}^{i+1,1})(\mE_{k-1}^{11}+\mC_{k-1}^{11})^{-1}\\[3mm]
\nonumber&&\cdot(\mE_{k-1}^{1,j+1}+\mC_{k-1}^{1,j+1}).
\end{eqnarray}

\item $l_3^\prime<l_2^\prime+1$\\
We decompose $\mX_k$ as
\begin{eqnarray}
 \nonumber\mX_k=(\mX_{k-l_2^\prime}^\prime,\ldots,\vx_{k-l_3^\prime+2}^\prime,\ldots ,\vx_{k-1}^\prime,\vx_k^\prime)^\prime
\end{eqnarray}
and $J(\mX_k)$ correspondingly as
\begin{flushleft}
$J(\mX_k)=$
\end{flushleft}
\begin{eqnarray}
\nonumber\left(
\begin{array}{ccccc}
\mA_k^{11} & \cdots &\mA_k^{1,l_{233}^\prime}&\cdots& \mA_k^{1,l_{21}^\prime}\\
\vdots &\ldots& \vdots &\cdots& \vdots \\
\mA_k^{l_{233}^\prime,1}&\cdots&\mA_k^{l_{233}^\prime,l_{233}^\prime}&\cdots&\mA_k^{l_{233}^\prime,l_{21}^\prime}\\
\vdots &\ldots& \vdots &\cdots& \vdots \\
\mA_k^{l_{21}^\prime,1} & \cdots & \mA_k^{l_{21}^\prime,l_{233}^\prime}&\cdots&\mA_k^{l_{21}^\prime,l_{21}^\prime}\\
\end{array}
\right)
\end{eqnarray}
where
\begin{eqnarray}\label{Eqs_secA_03}
\nonumber\mA_k^{1,1}&=&E(-\bigtriangleup_{\mX_{k-l_2^\prime}}^{\mX_{k-l_2^\prime}}\ln~p_k),\\
\nonumber\mA_k^{1,j}&=&E(-\bigtriangleup_{\mX_{k-l_2^\prime}}^{\vx_{k-l_2^\prime+j-1}}\ln~p_k),\\
\nonumber\mA_k^{i,1}&=&E(-\bigtriangleup_{\vx_{k-l_2^\prime+i-1}}^{\mX_{k-l_2^\prime}}\ln~p_k),\\
\nonumber\mA_k^{i,j}&=&E(-\bigtriangleup_{\vx_{k-l_2^\prime+i-1}}^{\vx_{k-l_2^\prime+j-1}}\ln~p_k),\\
 &&for~i,j=2,\ldots,l_2^\prime+1.
\end{eqnarray}

Using (\ref{Eqs_sec2_8}), (\ref{Eqsm_32}), (\ref{Eqsm_33}), the posterior information matrix for $\mX_{k+1}$ can be written in block form as
\begin{flushleft}
$J(\mX_{k+1})=$
\end{flushleft}
\begin{eqnarray}
\nonumber\left(
  \begin{array}{cccc}
    \mA_k^{1,1} & \ldots & \mA_k^{1,l_{21}^\prime} & 0 \\
    \ldots & \ldots & \ldots & \ldots \\
    \mA_k^{l_{21}^\prime,1} & \ldots & \mA_k^{l_{21}^\prime,l_{21}^\prime} & \mC_k^{l_{31}^\prime,l_3^\prime} \\
    ~ & ~ & +\mC_k^{l_{31}^\prime,l_{31}^\prime} & +\mB_k^{l_2^\prime,l_{21}^\prime} \\
    ~&~&+\mB_k^{{l_2^\prime},l_2^\prime}&~\\[3mm]
    0 & \cdots &\mC_k^{l_3^\prime,l_{31}^\prime}& \mC_k^{l_3^\prime,l_3^\prime} \\
~&~&+\mB_k^{l_{21}^\prime,l_2^\prime}&+\mB_k^{l_{21}^\prime,l_{21}^\prime}\\
  \end{array}
\right)
\end{eqnarray}
where $\mB_k^{i,j}$ and $\mC_k^{i,j}$ are defined in (\ref{Eqs_sec3_3})-(\ref{Eqs_sec3_4}). To simplify the matrix $J(\mX_{k+1})$, we denote by
$l_{233}^\prime\triangleq l_2^\prime-l_3^\prime+3$, $l_{232}^\prime\triangleq l_2^\prime-l_3^\prime+2$, $l_{31}^\prime\triangleq l_3^\prime-1$,
$l_{21}^\prime\triangleq l_2^\prime+1$. The information submatrix $\mJ_{k+1}$ for estimating $\vx_{k+1}\in R^r$ is given as the inverse of the
($r\times r$) right-lower block of $[J(\mX_{k+1})]^{-1}$, i.e.,
\begin{eqnarray}
\label{Eqs_secA_05}\mJ_{k+1}&=& \mD_k^{22}-\mD_k^{21}(\mE_k+\mD_k^{11})^{-1}\mD_k^{12},\\
\mE_k^{ij}&=&\mA_k^{i+1j+1}-\mA_k^{i+1,1}(\mA_k^{11})^{-1}\mA_k^{1,j+1}\\
\nonumber&=&(\mE_k^{ji})^\prime,\\
\nonumber &&for~i,j=1,\ldots,l_2^\prime,
\end{eqnarray}
Since  $\mD_k^{11}$, $\mD_k^{12}$, $\mD_k^{21}$ and $\mD_k^{22}$ in (\ref{Eqs_secA_05}) can be calculated as
(\ref{Eqs_sec3_5})-(\ref{Eqs_sec3_7}), we only needs to derive the recursion $\mE_k^{i,j}$.

Based on the  recursion between $p_k$ and $p_{k-1}$ given in (\ref{Eqs_sec2_8}) which depends on (\ref{Eqsm_32}) and (\ref{Eqsm_33}), and using
the definitions of $\mA_k^{i,j}$, $\mB_k^{i,j}$, $\mC_k^{i,j}$ given in (\ref{Eqs_secA_03}), (\ref{Eqs_sec3_3}), (\ref{Eqs_sec3_4}), respectively, it follows that
\begin{flushleft}
$ \left(\begin{array}{cc}
\mA_k^{11} & \mA_k^{1,j+1} \\
\mA_k^{i+1,1}& \mA_k^{i+1,j+1} \\
\end{array}\right) $
\end{flushleft}
\begin{eqnarray}\label{Eqsm_23}
=\left(
\begin{array}{ccc}
\mA_{k-1}^{11} & \mA_{k-1}^{12} & \mA_{k-1}^{1,j+2} \\[3mm]
\mA_{k-1}^{21} & \mA_{k-1}^{22}&\mA_{k-1}^{2,j+2} \\
~&+\mB_{k-1}^{11}&+\mB_{k-1}^{1,j+1}\\[3mm]
\mA_{k-1}^{i+2,1} & \mA_{k-1}^{i+2,2} &
\mA_{k-1}^{i+2,j+2}\\
~&+\mB_{k-1}^{i+1,1}&+\mB_{k-1}^{i+1,j+1}\\
~&~&+\mC_{k-1}^{i-l_2^\prime+l_3^\prime,j-l_2^\prime+l_3^\prime}\\
\end{array}
\right),
\end{eqnarray}
\begin{center}
$for~i=1, \ldots, l_2^\prime, ~j=1, \ldots, l_2^\prime,$
\end{center}
where if $i\leq l_2^\prime-l_3^\prime$ or $j\leq l_2^\prime-l_3^\prime$, then $\mC_{k-1}^{i-l_2^\prime+l_3^\prime,j-l_2^\prime+l_3^\prime}=0$;
if $i>l_2^\prime+1$ or $j>l_2^\prime+1$, then $\mA_{k-1}^{i,j}=0$.
Moreover, similar to the derivation of Equation (\ref{Eqs_secA_01}), we have the recursion of the matrix $\mE_k$ as
\begin{eqnarray}
\nonumber
\mE_k^{ij}&=&\mE_{k-1}^{i+1,j+1}+\mC_{k-1}^{i-l_2^\prime+l_3^\prime,j-l_2^\prime+l_3^\prime}+\mB_{k-1}^{i+1,j+1}\\[3mm]
\label{Eqs_secA_02}&&-(\mE_{k-1}^{i+1,1}+\mC_{k-1}^{i+1,1})(\mE_{k-1}^{11}+\mC_{k-1}^{11})^{-1}\\[3mm]
\nonumber&&\cdot(\mE_{k-1}^{1,j+1}+\mC_{k-1}^{1,j+1}).
\end{eqnarray}

\item $l_2^\prime+1=l_3^\prime$\\
We decompose $\mX_k$ as $\mX_k$=$(\mX_{k-l_3^\prime+1}^\prime,\ldots,\vx_{k-1}^\prime,\vx_k^\prime)^\prime$ and $J(\mX_k)$ correspondingly as
\begin{eqnarray}
\nonumber J(\mX_k)=\left(
\begin{array}{ccc}
\mA_k^{11} & \cdots & \mA_k^{1,l_3^\prime}\\
\vdots &\ldots& \vdots \\
\mA_k^{l_3^\prime,1} & \cdots &\mA_k^{l_3^\prime,l_3^\prime}\\
\end{array}
\right)
\end{eqnarray}
where
\begin{eqnarray}
\label{Eqs_secA_003}
\nonumber\mA_k^{1,1}&=&E(-\bigtriangleup_{\mX_{k-l_3^\prime+1}}^{\mX_{k-l_3^\prime+1}}\ln~p_k),\\
\nonumber\mA_k^{1,j}&=&E(-\bigtriangleup_{\mX_{k-l_3^\prime+1}}^{\vx_{k-l_3^\prime+j}}\ln~p_k),\\
\nonumber\mA_k^{i,1}&=&E(-\bigtriangleup_{\vx_{k-l_3^\prime+i}}^{\mX_{k-l_3^\prime+1}}\ln~p_k),\\
\nonumber\mA_k^{i,j}&=&E(-\bigtriangleup_{\vx_{k-l_3^\prime+i}}^{\vx_{k-l_3^\prime+j}}\ln~p_k),\\
 &&for ~i,j=2,\ldots,l_3^\prime.
\end{eqnarray}

Using (\ref{Eqs_sec2_8}), (\ref{Eqsm_32}), (\ref{Eqsm_33}), the posterior information matrix for $\mX_{k+1}$ can be written in block form as
\begin{flushleft}
$J(\mX_{k+1})=$
\end{flushleft}
\begin{eqnarray}
\nonumber \left(
\begin{array}{ccccc}
\mA_k^{1,1} &\mA_k^{1,2}&\cdots& \mA_k^{1,l_3^\prime} & 0 \\[3mm]

\mA_k^{2,1} & \mA_k^{2,2} &\cdots&\mA_k^{2,l_3^\prime} &\mC_k^{1,l_3^\prime} \\
~&+\mC_k^{1,1}&~&+\mC_k^{1,l_3^\prime-1}&+\mB_k^{1,l_3^\prime}\\
~&+\mB_k^{1,1}&~&+\mB_k^{1,l_3^\prime-1}&~\\
\vdots &\cdots & \vdots &\cdots& \vdots \\
\mA_k^{l_3^\prime,1}&\mA_k^{l_3^\prime,2}&\cdots&\mA_k^{l_3^\prime,l_3^\prime}&
\mC_k^{l_3^\prime-1,l_3^\prime}\\
~&+\mC_k^{l_3^\prime-1,1}&~&+\mC_k^{l_3^\prime-1,l_3^\prime-1}&+\mB_k^{l_3^\prime-1,l_3^\prime}\\
~&+\mB_k^{l_3^\prime-1,1}&~&+\mB_k^{l_3^\prime-1,l_3^\prime-1}&~\\[3mm]
0 & \mC_k^{l_3^\prime,1} &\cdots&
\mC_k^{l_3^\prime,l_3^\prime-1}& \mC_k^{l_3^\prime,l_3^\prime}\\
~&+\mB_k^{l_3^\prime,1}&~&+\mB_k^{l_3^\prime,l_3^\prime-1}&+\mB_k^{l_3^\prime,l_3^\prime} \\
\end{array}
\right)
\end{eqnarray}
where  $\mB_k^{i,j}$ and $\mC_k^{i,j}$ are defined in (\ref{Eqs_sec3_3})-(\ref{Eqs_sec3_4}). The information submatrix $\mJ_{k+1}$ for
estimating $\vx_{k+1}\in R^r$ is given as the inverse of the ($r\times r$) right-lower block of $[J(\mX_{k+1})]^{-1}$, i.e.,
\begin{eqnarray}
\label{Eqs_secA_07}\mJ_{k+1}&=& \mD_k^{22}-\mD_k^{21}(\mE_k+\mD_k^{11})^{-1}\mD_k^{12},\\
\mE_k^{ij}&=&\mA_k^{i+1j+1}-\mA_k^{i+1,1}(\mA_k^{11})^{-1}\mA_k^{1,j+1},\\
\nonumber&=&(\mE_k^{ji})^\prime,\\
\nonumber &&for ~i,j=1,\ldots,l_3^\prime-1.
\end{eqnarray}
Since  $\mD_k^{11}$, $\mD_k^{12}$, $\mD_k^{21}$ and $\mD_k^{22}$ in (\ref{Eqs_secA_07}) can be calculated as
(\ref{Eqs_sec3_12})-(\ref{Eqs_sec3_14}), we only needs to derive the recursion $\mE_k^{i,j}$.

Based on the recursion between $p_k$ and $p_{k-1}$ given in (\ref{Eqs_sec2_8}) which depends on (\ref{Eqsm_32}) and (\ref{Eqsm_33}), and using
the definitions of $\mA_k^{i,j}$, $\mB_k^{i,j}$, $\mC_k^{i,j}$ given in (\ref{Eqs_secA_003}), (\ref{Eqs_sec3_3}), (\ref{Eqs_sec3_4}), respectively, it follows that
\begin{flushleft}
$\left(
\begin{array}{cc}
\mA_k^{11} & \mA_k^{1,j+1} \\
\mA_k^{i+1,1}& \mA_k^{i+1,j+1} \\
\end{array}
\right)$
\end{flushleft}
\begin{eqnarray}\label{Eqsm_28}
\nonumber =\left(
\begin{array}{ccc}
\mA_{k-1}^{11} & \mA_{k-1}^{12} & \mA_{k-1}^{1,j+2} \\[3mm]
\mA_{k-1}^{21} & \mA_{k-1}^{22}&\mA_{k-1}^{2,j+2} \\
~&+\mC_{k-1}^{11}&+\mC_{k-1}^{1,j+1}\\
~&+\mB_{k-1}^{11}&+\mB_{k-1}^{1,j+1}\\[3mm]
\mA_{k-1}^{i+2,1} & \mA_{k-1}^{i+2,2} & \mA_{k-1}^{i+2,j+2}\\
~&+\mC_{k-1}^{i+1,1}&+\mC_{k-1}^{i+1,j+1}\\
~&+\mB_{k-1}^{i+1,1}&+\mB_{k-1}^{i+1,j+1}
\end{array}
\right),
\end{eqnarray}
\begin{center}
$ for ~i=1, \ldots, l_3^\prime-1, j=1, \ldots, l_3^\prime-1$.
\end{center}
where $\mA_{k-1}^{l_3^\prime+1,j}=0$ and $\mA_{k-1}^{i,l_3^\prime+1}=0$. Moreover, similar to the derivation of Equation (\ref{Eqs_secA_01}), we
have the recursion of the matrix $\mE_k$ as

\begin{eqnarray}
\nonumber
\mE_k^{ij}&=&\mE_{k-1}^{i+1,j+1}+\mC_{k-1}^{i+1,j+1}+\mB_{k-1}^{i+1,j+1}\\[3mm]
\nonumber&&-(\mE_{k-1}^{i+1,1}+\mB_{k-1}^{i+1,1}+\mC_{k-1}^{i+1,1})\\[3mm]
\label{Eqs_secA_04}&&\cdot(\mE_{k-1}^{11}+\mB_{k-1}^{11}+\mC_{k-1}^{11})^{-1}\\[3mm]
\nonumber&&\cdot(\mE_{k-1}^{1,j+1}+\mB_{k-1}^{1,j+1}+\mC_{k-1}^{1,j+1}).
\end{eqnarray}
\end{enumerate}
Based on (\ref{Eqs_secA_06}), (\ref{Eqs_secA_01}), (\ref{Eqs_secA_05}), (\ref{Eqs_secA_02}), (\ref{Eqs_secA_07}) and (\ref{Eqs_secA_04}), we
have completed the proof of the Theorem \ref{thm_1}.

\subsection{proof of Corollary \ref{coro_2}}
In case of $l_1=0$, $l_2=0$, $l_4=0$, $l_3=l$,  by (\ref{Eqsm_32})-(\ref{Eqsm_33}), Equation (\ref{Eqs_sec2_8}) can be written as
\begin{flushleft}
$p_{k+1}=$
\end{flushleft}
\begin{eqnarray}\label{Eqsm_34}
\nonumber\left\{
\begin{array}{ll}
p_kp(\vx_{k+1}|\vx_k)p(\vz_{k+1}|\vx_{k+1}), & \hbox{for ~$l=1$ } \\
~&\hbox{or $l=0$,}\\
p_kp(\vx_{k+1}|\vx_k)p(\vz_{k+1}|\vx_{k+1},\vx_{k})& \hbox{for ~ $l=2$,}\\
p_kp(\vx_{k+1}|\vx_k)p(\vz_{k+1}|\vx_{k+1},\ldots,\vx_{k-l+2})& \hbox{ for ~ $l>2$}\\
\end{array}
\right.
\end{eqnarray}
Thus, we can immediately obtain the recursion (\ref{Eqs_sec3_15}) by Theorem \ref{thm_1}. The recursion of $\mE_k$ can be written as the three
cases of (\ref{Eqs_sec3_16}), (\ref{Eqs_sec3_17}) and (\ref{Eqs_sec3_18}), respectively. At the same time, the matrices $\mB_k$, $\mC_k$, and
$\mD_k$ become correspondingly appropriate forms.

\subsection{proof of Corollary \ref{coro_3}}
In case of $l_1=0$, $l_2=l$, $l_3=0$, $l_4=0$, i.e., $l_3^\prime<l_2^\prime+1$, by (\ref{Eqsm_32})-(\ref{Eqsm_33}), Equation (\ref{Eqs_sec2_8})
can be  written as
\begin{eqnarray}
\nonumber p_{k+1}=p_kp(\vx_{k+1}|\vx_k,\ldots,\vx_{k-l_2^\prime+1})p(\vz_{k+1}|\vx_{k+1}).
\end{eqnarray}
Thus, we can get (\ref{Eqs_sec3_19})-(\ref{Eqs_sec3_20}) by Theorem \ref{thm_1}.  At the same time, the matrices $\mB_k$, $\mC_k$, and $\mD_k$
become correspondingly appropriate forms.

\subsection{proof of Corollary \ref{coro_4}}
In case of $l_1=l$, $l_2=0$, $l_3=0$, $l_4=0$, i.e., $l_3^\prime<l_2^\prime+1$, by (\ref{Eqsm_33}), Equation (\ref{Eqs_sec2_8}) can be
simplified as
\begin{eqnarray}
\nonumber p_{k+1}=\left\{
\begin{array}{ll}
p_kp(\vx_{k+1}|\vx_k)p(\vz_{k+1}|\vx_{k+1})& \hbox{for ~ $l=0$},\\
p_kp(\vx_{k+1}|\vx_k)p(\vz_{k+1}|\vx_{k+1},\\
\qquad\qquad\qquad\vz_k,\ldots,\vz_{k-l+1})& \hbox{ for ~$l>0$}.\\
\end{array}
\right.
\end{eqnarray}
Thus, we can get (\ref{Eqs_sec3_24}) by Theorem \ref{thm_1}. At the same time, the matrices $\mB_k$, $\mC_k$, and $\mD_k$ become correspondingly
appropriate forms.


\begin{thebibliography}{10}

\bibitem{Li-Jilkov03}
X.~R. Li and V.~P. Jilkov, ``Survey of maneuvering target tracking. {Part I:}
  dynamic models,'' {\em IEEE Transactions on Aerospace and Electronic
  Systems}, vol.~39, no.~4, pp.~1333--11364, 2003.

\bibitem{Mazor-Averbuch-BarShalom-Dayan98}
E.~Mazor, A.~Averbuch, Y.~Bar-Shalom, and J.~Dayan, ``Interacting multiple
  model methods in target tracking: {A} survey,'' {\em IEEE Transactions on
  Aerospace and Electronic Systems}, vol.~34, pp.~103--123, JANUARY 1998.

\bibitem{Simon06}
D.~Simon, {\em Optimal State Estimation: {Kalman}, $H_\infty$, and Nonlinear
  Approaches}.
\newblock Wiley-Interscience, 2006.

\bibitem{BarShalom-Li-Kirubarajan01}
Y.~Bar-Shalom, X.~Li, and T.~Kirubarajan, {\em Estimation with Applications to
  Tracking and Navigation}.
\newblock New York: Wiley, 2001.

\bibitem{Ljung-Gunnarsson90}
L.~Ljung and S.~Gunnarsson, ``Adaptation and tracking in system
  identification--a survey,'' {\em Automatica}, vol.~26, pp.~7--21, 1990.

\bibitem{Guo-94}
L.~Guo, ``Stability of recursive stochastic tracking algorithms,'' {\em SIAM
  Journal on Control and Optimization}, vol.~32, pp.~1195--1225, 1994.

\bibitem{Maryak-Spall-Silberman-95}
J.~L. Maryak, J.~C. Spall, and G.~L. Silberman, ``Uncertainties for recursive
  estimators in nonlinear state-space models, with applications to
  epidemiology,'' {\em Automatica}, vol.~31, no.~12, pp.~1889--1892, 1995.

\bibitem{Rogers87}
S.~R. Rogers, ``Alpha-beta filter with correlated measurement noise,'' {\em
  IEEE Transactions on Aerospace and Electronic Systems}, vol.~23,
  pp.~592--594, July 1987.

\bibitem{Halevi90}
Y.~Halevi, ``Optimal observers for systems with colored noises,'' {\em IEEE
  Transactions on Automatic Control}, vol.~35, pp.~1075--1078, August 1990.

\bibitem{Blair-Watson-Rice91}
W.~D. Blair, G.~A. Watson, and T.~R. Rice, ``Tracking maneuvering targets with
  an interacting multiple model filter containing exponentially correlated
  acceleration models,'' in {\em Proceedings of the Twenty-Third Southeastern
  Symposium on System Theory}, pp.~224--228, 1991.

\bibitem{Rapoport-Oshman05}
I.~Rapoport and Y.~Oshman, ``A {Cram$\acute{e}$r-Rao-Type} estimation lower
  bound for systems with measurement faults,'' {\em IEEE Transactions on
  Automatic Control}, vol.~50, pp.~1234--1245, September 2003.

\bibitem{Yun-Bachmann06}
X.~Yun and E.~R. Bachmann, ``Design, implementation, and experimental results
  of a quaternion-based {Kalman} filter for human body motion tracking,'' {\em
  IEEE Transactions on Robtics}, vol.~22, pp.~1216--1227, December 2006.

\bibitem{Jiang-Zhou-Zhu10}
P.~Jiang, J.~Zhou, and Y.~M. Zhu, ``Globally optimal {Kalman} filtering with
  finite-time correlated noises,'' in {\em Proceedings of the 49th IEEE
  Conference on Decision and Control}, pp.~15--17, December 2010.

\bibitem{Chen12}
S.~Y. Chen, ``Kalman filter for robot vision: {A} survey,'' {\em IEEE
  Transactions on Industrial Electronics}, vol.~59, pp.~4409--4420, November
  2012.

\bibitem{VanTrees68}
H.~L.{Van Trees}, {\em Detection, Estimation, and Modulation Theory, Part I.}
\newblock New York: Wiley, 1968.

\bibitem{BarShalom2014}
Bar-Shalom, Y.~Osborne, R.~Willett, and F.~P.~Daum, ``{CRLB} for likelihood
  functions with parameter-dependent support and a new bound,'' {\em Aerospace
  Conference, 2014 IEEE}, March 2014.

\bibitem{Tichavsky-Muravchik-Nehorai98}
P.~Tichavsk$\acute{y}$, C.~H. Muravchik, and A.~Nehorai, ``Posterior
  {Cram$\acute{e}$r-Rao} bounds for discrete-time nonlinear filtering,'' {\em
  IEEE Transactions on Signal Processing}, vol.~46, pp.~1386--1396, May 1998.

\bibitem{Weinstein-Weiss88}
E.~Weinstein and A.~J. Weiss, ``A general class of lower bounds in parameter
  estimation,'' {\em IEEE Transactions on Information Theory}, vol.~34,
  pp.~338--342, March 1988.

\bibitem{Bobrovsky-Zakai75}
B.~Z. Bobbovsky and M.~Zakai, ``A lower bound on the estimation error for
  {Markov} processes,'' {\em IEEE Transactions on Automatic Control}, vol.~20,
  pp.~785--788, 1975.

\bibitem{Taylor79}
J.~H. Taylor, ``The {Cram$\acute{e}$r-Rao} estimation error lower bound
  computation for deterministic nonlinear systems,'' {\em IEEE Transactions on
  Automatic Control}, vol.~24, pp.~343--344, April 1979.

\bibitem{Galdos80}
J.~I. Galdos, ``A {Cram$\acute{e}$r-Rao} bound for multidimensional
  discrete-time dynamical systems,'' {\em IEEE Transactions on Automatic
  Control}, vol.~25, pp.~117--119, 1980.

\bibitem{Chang81}
C.~B. Chang, ``Two lower bounds on the covariance for nonlinear estimation
  problems,'' {\em IEEE Transactions on Automatic Control}, vol.~26,
  pp.~1294--1297, December 1981.

\bibitem{Kerr89}
T.~H. Kerr, ``Status of {CR-like} lower bounds for nonlinear filtering,'' {\em
  IEEE Transactions on Aerospace and Electronic Systems}, vol.~25,
  pp.~590--600, September 1989.

\bibitem{Koshaev-Stepanov97}
D.~A. Koshaev and O.~A. Stepanov, ``Application of the {Rao-Cramer} inequality
  in problems of nonlinear estimation,'' {\em Computer and Systems Sciences
  International}, vol.~36, no.~2, pp.~220--227, 1997.

\bibitem{Schultheiss-Weinstein81}
P.~M. Schultheiss and E.~Weinstein, ``Lower bounds on the localization errors
  of a moving source observed by a passive array,'' {\em IEEE Transactions on
  Acoustics, Speech, and Signal Processing}, vol.~29, pp.~600--607, June 1981.

\bibitem{Aidala-Hammel83}
V.~J. Aidala and S.~E. Hammel, ``Utilization of modified polar coordinates for
  bearings-only tracking,'' {\em IEEE Transaction on Automatic Control},
  vol.~28, pp.~283--294, March 1983.

\bibitem{Kirubarajan-BarShalom96}
T.~Kirubarajan and Y.~Bar-Shalom, ``Low observable target motion analysis using
  amplitude information,'' {\em IEEE Transactions on Aerospace and Electronic
  Systems}, vol.~32, pp.~1367--1384, October 1996.

\bibitem{Niu-Willett-BarShalom01}
R.~Niu, P.~Willett, and Y.Bar-shalom, ``Matrix {CRLB} scaling due to
  measurements of uncertain origin,'' {\em IEEE Transactions on Signal
  Processing}, vol.~49, pp.~1325--1335, July 2001.

\bibitem{Zhang-Willett-BarShalom02}
X.~Zhang, P.~Willett, and Y.~Bar-Shalom, ``The {Cram$\acute{e}$r-Rao} bound for
  dynamic target tracking with measurement origin uncertainty,'' in {\em
  Proceedings of the 41st IEEE Conference on Decision and Control}, (Las
  Vegas), pp.~3428--3433, December 2002.

\bibitem{Zheng-Ozdemir-Niu-Varshney12}
Y.~Zheng, O.~Ozdemir, R.~Niu, and P.~K. Varshney, ``New conditional posterior
  {Cramer-Rao} lower bounds for nonlinear sequential {Bayesian} estimation,''
  {\em IEEE Transactions on Signal Processing}, vol.~60, pp.~5549--5556,
  October 2012.

\bibitem{Kar-Varshney-Palaniswami12}
S.~Kar, P.~K. Varshney, and M.~Palaniswami, ``{Cramer-Rao} bounds for
  polynomial signal estimation using sensors with {AR(1)} drift,'' {\em IEEE
  Transactions on Signal Processing}, vol.~60, pp.~5494--5507, October 2012.

\bibitem{Joshi-Boyd09}
S.~Joshi and S.~Boyd, ``Sensor selection via convex optimization,'' {\em IEEE
  Transactions on Signal Processing}, vol.~57, pp.~451--462, February 2009.

\bibitem{Shen-Varshney14}
X.~Shen and P.~K. Varshney, ``Sensor selection based on generalized information
  gain for target tracking in large sensor networks,'' {\em IEEE Transactions
  on Signal Processing}, vol.~62, pp.~363--375, January 2014.

\bibitem{BarShalom02}
Y.~Bar-Shalom, ``Update with out-of-sequence measurements in tracking: Exact
  solution,'' {\em IEEE Transactions on Aerospace and Electronic Systems},
  vol.~38, no.~3, pp.~769--778, 2002.

\bibitem{Horn-Johnson12}
R.~A. Horn and C.~R. Johnson, {\em Matrix Analysis}.
\newblock Cambridge University Press, 2nd revised~ed., 2012.

\end{thebibliography}
\end{document}